\documentclass[twoside, a4paper, leqno]{article}

\usepackage{amsmath, amssymb, amsthm}
\usepackage{mathtools}
\usepackage{mathrsfs}
\usepackage[sort&compress, square, comma, numbers]{natbib} 
\usepackage{hyperref}
\usepackage{bm}
\usepackage{graphicx}        
\usepackage{multicol}        
\usepackage[bottom]{footmisc}
\usepackage{comment}
\usepackage{geometry}

\newcommand{\D}{\mathbf{D}}

\newcommand{\norm}[1]{\left\lVert{\textstyle #1 }\right\rVert}

\graphicspath{   
{FIGURES/}
}

\begin{document}

\title{Computational design of nanophotonic structures using an
adaptive  finite element method}

\author{L. Beilina \thanks{ Department of Mathematical Sciences,
    Chalmers University of Technology and Gothenburg University,
    SE-42196 Gothenburg, Sweden, e-mail:
    \texttt{\ larisa@chalmers.se}} \and L. Mpinganzima \thanks{
    Department of Physics, Chalmers University of Technology, SE-42196
    Gothenburg, Sweden, e-mail:
    \texttt{\ lydie.mpinganzima@chalmers.se}} \and P. Tassin \thanks{
    Department of Physics, Chalmers University of Technology and
    Gothenburg University, SE-42196 Gothenburg, Sweden, e-mail:
    \texttt{\ nanophotonicstructures@gmail.com}}}

\date{}

\maketitle

\begin{abstract}

We consider the problem of the construction of the nanophotonic
structures of arbitrary geometry with prescribed desired properties.
We reformulate this problem as an optimization problem for the
Tikhonov functional which is minimized on adaptively locally refined
meshes. These meshes are refined only in places where the
nanophotonic structure should be designed. Our special symmetric mesh
refinement procedure allows the construction of different nanophotonic
structures.
We illustrate efficiency of our adaptive optimization algorithm on the
construction of nanophotonic structure  in two dimensions.

\end{abstract}



\section{Introduction}

The goal of this work is to develop a new optimization algorithm that
can construct arbitrary nanophotonic structures from desired
scattering parameters.  Nanophotonics is the study of the interaction
of electromagnetic waves with structures that have feature sizes equal
or smaller than the wavelength of the waves. Examples are photonic
crystals (structured on the wavelength scale), metamaterials
(subwavelength structured media with new optical properties that are
not available from natural materials) and plasmonic
devices (exploiting collective excitations in metals that result in
strong field enhancement) \cite{JoanJohnson, Maier,  SoukoulisW, ZheludevK}.

 In this paper, we present a nonparametric optimization algorithm that
 can find inner structure of the domain with arbitrary geometry. To do
 that we apply an adaptive finite element method of \cite{B,BJ} with
 iterative choice of the regularization parameter  \cite{BKS}.  We
 illustrate the efficiency of the proposed adaptive optimization
 method on the solution of the hyperbolic coefficient inverse problem
 (CIP) in two dimensions.  The goal of our numerical simulations is to
 reconstruct the permittivity function of the hyperbolic equation from
 single observations of the transmitted and backscattered solution of
 this equation in space and time.  For computational solution of this
 inverse problem we use the domain decomposition method of
 \cite{hybrid}. To solve our CIP we minimize the corresponding
 Tikhonov functional via Lagrangian approach.  This approach is
 similar to the one applied recently in \cite{B,  HybEf, BJ, BCN, BN} for
 the solution of different hyperbolic CIPs: we find optimality
 conditions which express stationarity of the Lagrangian, involving
 the solution of state and adjoint equations together with an equation
 expressing that the gradient of the Lagrangian with respect to the
 permittivity function vanishes. Then we construct an adaptive
 conjugate gradient algorithm and compute the unknown permittivity
 function in an iterative process by solving in every step the state
 and adjoint hyperbolic equations and updating in this way the desired
 permittivity function.

\section{Statement of the forward and inverse problems}

\label{sec:modelhyb}

 Let $x = (x_1, x_2)$ denotes a point in $\mathbb{R}^2$ in an
 unbounded domain $D$.  In this work we consider the propagation of
 electromagnetic waves in two dimensions with a field
 polarization. Thus, we model the wave propagation by the following
 Cauchy problem for the scalar wave equation:
\begin{equation}\label{modelhyb}
\begin{cases}
\varepsilon(x) \frac{\partial^2 E}{\partial t^2}  -  \triangle E  = \delta(x_2 - x_0)p(t) & ~ \mbox{in}~~\mathbb{R}^2 \times (0, \infty), \\
  E(x,0) = f_0(x), ~~~E_t(x,0) = 0 &~ \mbox{in}~~ D.
\end{cases}
\end{equation}
Here, $E$ is the electric field generated by the plane wave $p(t)$
which is incident at $x_2 = x_0$ and propagates along $x_2$ axis,
$\varepsilon(x)$ is the spatially distributed dielectric
permittivity. We note that in this work we use the single equation
(\ref{modelhyb}) instead of the full Maxwell's equations, since in
\cite{BMaxwell} was demonstrated numerically that in the similar
numerical setting, as we will use in this note, other components of
the electric field are negligible compared to the initialized one.  We
also note that a scalar model of the wave equation was used
successfully to validate reconstruction of the dielectric permittivity
function with transmitted \cite{BK1, BK2} and backscattered
experimental data \cite{ btkm14, btkm14b, KBKSNF, NBKF, NBKF2}.

Let now $D \subset \mathbb{R}^{2}$ be a convex bounded domain
with the boundary $\partial D \in C^{2}$.  We denote by
$D_T := D \times (0,T), \partial D_T := \partial D \times (0,T), T > 0$ and assume that
\begin{equation}
 f_{0}\in H^{1}(D),  \varepsilon(x) \in C^2(D).   \label{f1} 
\end{equation}
For computational solution of (\ref{modelhyb}) we use the domain
decomposition finite element/finite difference (FE/FD) method of
\cite{hybrid, BSA} which was applied for the solution of different
coefficient inverse problems for the acoustic wave equation in
\cite{B, hybrid, BJ}.  To apply method of \cite{hybrid, BSA} we
decompose $D$ into two regions $D_{FEM}$ and $D_{FDM}$ such that the
whole domain $D = D_{FEM} \cup D_{FDM}$, and $D_{FEM} \cap D_{FDM} =
\emptyset$.  In $D_{FEM}$ we use the finite element method (FEM), and
in $D_{FDM}$ we will use the Finite Difference Method (FDM).  We avoid
instabilities at interfaces between FE and FD domains since FE and FD
discretization schemes coincide on two common structured layers with
$\varepsilon(x)=1$ in them.

 Let the boundary $\partial D$ be such that $\partial D =\partial _{1}
 D \cup \partial _{2} D \cup \partial _{3} D$ where $\partial _{1} D$
 and $\partial _{2} D$ are, respectively, front and back sides of the
 domain $D$, and $\partial _{3} D$ is the union of left, right, top
 and bottom sides of this domain.  At $S_{T_1} := \partial_1 D
 \times (0,T)$ and $S_{T_2} := \partial_2 D \times (0,T)$ we have
 time-dependent backscattering and transmission observations,
 correspondingly.  We  define $S_{1,1} := \partial_1 D \times
 (0,t_1]$, $S_{1,2} := \partial_1 D \times (t_1,T)$, and $S_3 :=
   \partial_3 D \times (0, T)$.
We
 also introduce the following spaces of real valued  functions
\begin{equation}\label{spaces}
\begin{split}
H_E^1(D_T) &:= \{ w \in H^1(D_T):  w( \cdot , 0) = 0 \}, \\
H_{\lambda}^1(D_T) &:= \{ w \in  H^1(D_T):  w( \cdot , T) = 0\},\\
U^{1} &=H_{E}^{1}(D_T)\times H_{\lambda }^{1}(D_T)\times C\left( \overline{D}\right),\\
\end{split}
\end{equation}
and define standard $L_2$ inner product and space-time norms, correspondingly, as
\begin{equation*}
\begin{split}
((u,v))_{D_T}    &= \int_{D} \int_0^T u v~ dx dt,~ \|u \|_{L_2(D_T)}^2  = ((u,u))_{D_T}, \\
(u,v)_{D}    &= \int_{D} u v ~dx,~ \| u \|_{L_2(D)}^2  = (u,u)_{D}.
\end{split}
\end{equation*}

In our computations we have used the following model problem 
\begin{equation}\label{model1}
\begin{cases}
\varepsilon \frac{\partial^2 E}{\partial t^2}  -  \triangle E   = 0 &~ \mbox{in}~~ D_T, \\
  E(x,0) = f_0(x), ~~~E_t(x,0) = 0 &~ \mbox{in}~~ D,     \\
\partial _{n} E = p\left( t\right) & ~\mbox{on}~ S_{1,1},
\\
\partial _{n}  E =-\partial _{t} E & ~\mbox{on}~ S_{1,2},
\\
\partial _{n} E =-\partial _{t} E &~\mbox{on}~ S_{T_2}, \\
\partial _{n} E =0 &~\mbox{on}~ S_3.\\
\end{cases}
\end{equation}
In (\ref{modelhyb}) we use the first order absorbing boundary
conditions \cite{EM}. These conditions are exact in the case of  computations
of section \ref{sec:Numer-Simul} since we initialize the plane wave orthogonal
to the domain of propagation.

We choose the coefficient $\varepsilon(x)$ in (\ref{modelhyb})  such that
\begin{equation} \label{coefic}
\begin{cases}
\varepsilon \left( x\right) \in \left ( 0, M \right], M=const. > 0,
& \text{ for }x\in  D _{FEM}, \\
 \varepsilon(x) =1
& \text{ for }x\in  D _{FDM}.
\end{cases}
\end{equation}

We consider the following inverse problem

\textbf{Inverse Problem (IP)} 

\emph{Let the coefficient }$\varepsilon \left( x\right)$\emph{\ in the problem (\ref{model1}) satisfy
conditions  (\ref{coefic}) and  assume that  }$ \varepsilon\left( x\right)
$\emph{\ is unknown in the domain }$D \diagdown
D_{FDM}$\emph{. Determine the function }$ \varepsilon\left( x\right)
$\emph{\  in (\ref{model1}) for }$x\in D \diagdown D_{FDM},$ \emph{\ assuming
  that the following function }$\widetilde E\left( x,t\right) $\emph{\ is
  known}
\begin{equation}
  E \left( x,t\right) = \widetilde E \left( x,t\right), ~\forall \left( x,t\right) \in S_{T_1} \cup S_{T_2}.  \label{2.5}
\end{equation}


\section{Optimization method}

\label{sec:opt}

In this section we present the reconstruction method to solve inverse
problem \textbf{IP}. This method is based on the finding of the
stationary point of the following Tikhonov functional 
\begin{equation}
F(E, \varepsilon) = \frac{1}{2} \int_{S_{T_1} \cup S_{T_2}}( E- \widetilde{E})^2 z_{\delta }(t) d\sigma dt +
\frac{1}{2} \gamma \int_{D}(\varepsilon- \varepsilon_g)^2~~ dx,
\label{functional}
\end{equation}
where $E$ satisfies the equations
(\ref{model1}), $\varepsilon_{0}$ is the initial guess
for $\varepsilon$,  $\widetilde{E}$ is the observed field at $S_{T_1} \cup S_{T_2}$, $\gamma > 0$ is the  regularization parameter and  $z_{\delta }$ can be
chosen  as in \cite{BCN}.

To find minimum of (\ref{functional}) we use the Lagrangian approach \cite{B, BJ} and define  the following Lagrangian
\begin{equation}\label{lagrangian1}
\begin{split}
L(v) &= F(E, \varepsilon) 
-  \int_{D_T} \varepsilon \frac{\partial
 \lambda }{\partial t} \frac{\partial E}{\partial t}  ~dxdt  
+   \int_{D_T}( \nabla  E)( \nabla  \lambda)~dxdt  \\
& - \int_{S_{1,1}} \lambda p(t) ~d \sigma dt    + \int_{S_{1,2}} \lambda \partial_t E ~d\sigma dt
   + \int_{S_{T_2}} \lambda \partial_t E ~d\sigma dt   , \\
\end{split}
\end{equation}
where $v=(E,\lambda, \varepsilon) \in U^1$, and search for a stationary point
with respect to $v$ satisfying $ \forall \bar{v}= (\bar{E}, \bar{\lambda}, \bar{\varepsilon}) \in U^1$
\begin{equation}
 L'(v; \bar{v}) = 0 ,  \label{scalar_lagr1}
\end{equation}
where $ L^\prime (v;\cdot )$ is the Jacobian of $L$ at $v$.

Similarly with \cite{B, BJ} we use conditions $\lambda \left( x,T\right) =\partial _{t}\lambda \left(
x,T\right) =0$ and imply such conditions on the function
$\lambda $  that $ L\left( E,\lambda, \varepsilon \right)
:=L\left( v\right) =F\left( E, \varepsilon\right).$ 
We also  use  conditions (\ref{coefic})  on $\partial D$, together with initial and boundary conditions  of (\ref{model1}) to get that
for all $\bar{v} \in U^1$,
\begin{equation}\label{forward1}
\begin{split}
0 &= \frac{\partial L}{\partial \lambda}(v)(\bar{\lambda}) =
- \int_{D_T} \varepsilon \frac{\partial \bar{\lambda}}{\partial t} \frac{\partial E}{\partial t}~ dxdt 
+  \int_{D_T}  (  \nabla E) (\nabla  \bar{\lambda}) ~ dxdt  \\
&- \int_{S_{1,1}} \bar{\lambda} p(t) ~d \sigma dt   + \int_{S_{1,2}} \bar{\lambda} \partial_t E ~d\sigma dt 
   + \int_{S_{T_2}} \bar{\lambda} \partial_t E  ~d\sigma dt,~~\forall \bar{\lambda} \in H_{\lambda}^1(D_T),
\end{split}
\end{equation}
\begin{equation} \label{control1}
\begin{split}
0 &= \frac{\partial L}{\partial E}(v)(\bar{E}) =
\int_{S_T}( E- \widetilde{E})~ \bar{E}~ z_{\delta}~ d \sigma dt- \int_{D} 
\varepsilon \frac{\partial{\lambda}}{\partial t}(x,0) \bar{E}(x,0) ~dx 
- \int_{S_{1,2} \cup S_{T_2}} 
\frac{\partial{\lambda}}{\partial t} \bar{E} ~d\sigma dt 
 \\
&-  \int_{D_T} \varepsilon  \frac{\partial \lambda}{\partial t} \frac{\partial \bar{E}}{\partial t}~ dxdt
 + \int_{D_T} ( \nabla  \lambda) (\nabla  \bar{E})  ~ dxdt, ~\forall \bar{E} \in H_{E}^1(D_T),
\end{split}
\end{equation}
\begin{equation} \label{grad1new} 
0 = \frac{\partial L}{\partial  \varepsilon}(v)(\bar{\varepsilon})
 =  -  \int_{D_T}  \frac{\partial \lambda}{\partial t} \frac{\partial E}{\partial t} \bar{\varepsilon}~dxdt 
+\gamma \int_{D} (\varepsilon - \varepsilon_g) \bar{\varepsilon}~dx,~ x \in D.
\end{equation}
We observe that (\ref{forward1}) is the weak formulation of the state equation
(\ref{model1}) and  (\ref{control1}) is the weak
formulation of the following adjoint problem
\begin{equation}\label{adjoint1}
\begin{cases} 
\varepsilon \frac{\partial^2 \lambda}{\partial t^2} - 
  \triangle \lambda  = -  (E - \widetilde{E}) z_{\delta} &~  x \in S_T,   \\
\lambda(\cdot, T) =  \frac{\partial \lambda}{\partial t}(\cdot, T) = 0, \\
\partial _{n} \lambda = \partial _{t} \lambda & ~\mbox{on}~ S_{1,2},
\\
\partial _{n} \lambda =\partial _{t} \lambda & ~\mbox{on}~ S_{T_2}, \\
\partial _{n} \lambda =0 & ~\mbox{on}~ S_3.
\end{cases}
\end{equation}

\section{Discretization of the domain decomposition FE/FD method}

\label{sec:fem}

As was mentioned above for the numerical solution of (\ref{modelhyb})
we use the domain decomposition FE/FD method of \cite{hybrid,
  BSA}. Similarly with these works, in our computations we decompose
the finite difference domain $D_{FDM}$ into squares and the finite
element domain $D_{FDM}$ - into triangles. For FDM discretization we
use the standard difference discretization of the equation
(\ref{model1}) and obtain an explicit scheme as in \cite{hybrid}.

For the finite element discretization of $D_{FEM}$ we define a
partition $K_{h}=\{K\}$ which consists of triangles. We define by $h$
the mesh function as $h|_{K}=h_{K}$, where $h_K$ is the local diameter of
the element $K$, and assume  the minimal
angle condition on the $K_{h}$ \cite{Brenner}. Let $J_{\tau }=\left\{ J\right\} $ be a partition of
the time interval $(0,\,T)$ into subintervals $J=(t_{k-1},\,t_{k}]$
of uniform length $\tau =t_{k}-t_{k-1}$.

To solve the state problem (\ref{model1}) and the
adjoint problem (\ref{adjoint1}) we define the finite element spaces, 
$W_{h}^{E}\subset H_{E}^{1}\left( D_{T}\right) $ and $W_{h}^{\lambda
}\subset H_{\lambda }^{1}\left( D_{T}\right) $. First, we introduce the
finite element trial space $W_{h}^{u}$ 
\begin{equation}
W_{h}^{u}:=\{w\in H_{u}^{1}(D_T):w|_{K\times J}\in P_{1}(K)\times
P_{1}(J),~\forall K\in K_{h},~\forall J\in J_{\tau }\},  
\end{equation}
where $P_{1}(K)$ and $P_{1}(J)$ denote the set of linear functions on
$K$ and $J$, respectively. We also introduce the finite element test
space $W_{h}^{\lambda }$  as
\begin{equation}
W_{h}^{\lambda }:=\{w\in H_{\lambda }^{1}(D_T):w|_{K\times J}\in
P_{1}(K)\times P_{1}(J),~\forall K\in K_{h},~\forall J\in J_{\tau }\}. 
\end{equation}%
To approximate the function $\varepsilon _\mathrm{r}$, we use the space of
piecewise constant functions $C_{h}\subset L_{2}\left( D \right) $, 
\begin{equation}
C_{h}:=\{u\in L_{2}(D ):u|_{K}\in P_{0}(K),~\forall K\in K_{h}\}, 
\end{equation}%
where $P_{0}(K)$ is the set of constant functions on $K$.

Setting $U_{h}=W_{h}^{E}\times W_{h}^{\lambda }\times C_{h}$, the
finite element method for  (\ref{scalar_lagr1}) now
reads: \emph{Find }$u_{h}\in U_{h}$\emph{, such that}
\begin{equation}
L^{\prime }(E_{h})(\bar{E})=0, ~\forall \bar{E}\in U_{h}.  
\end{equation}

\section{Adaptive conjugate gradient   algorithm}
\label{subsec:ad_alg}

To compute the minimum of the functional (\ref{functional}) we use the
adaptive conjugate gradient method (ACGM).  The regularization
parameter $\gamma$ in ACGM is computed iteratively via rules of
\cite{BKS}.  For the local mesh refinement we use a posteriori error
estimate  of \cite {B, BJ} which means that
the finite element mesh in $D_{FEM}$ should be locally refined where
the maximum norm of the Fr\'{e}chet derivative of the Lagrangian with
respect to the coefficient is large.

We denote
\begin{equation}\label{Bhm}
\begin{split}
  {g}^m(x) = - {\int_0}^T  \frac{\partial \lambda_h^m}{\partial t}
\frac{\partial E_h^m}{\partial t}~ dt   + \gamma^m (\varepsilon_h^m - \varepsilon_g), 
\end{split}
\end{equation}
where $\varepsilon_{h}^{m}$ is approximation of the function
$\varepsilon_{h}$ on the iteration $m$, $E_{h}\left(
x,t,\varepsilon_{h}^{m}\right) ,\lambda _{h}\left(
x,t,\varepsilon_{h}^{m} \right) $\ are computed by solving the state
(\ref{model1}) and adjoint (\ref{adjoint1}) problems, respectively,
with $\varepsilon:=\varepsilon_{h}^{m}$.

\textbf{Algorithm}

\begin{itemize}
\item Step 0.  Choose initial mesh $K_{h}$ in $D_{FEM}$ and time
  partition $J_{\tau}$ of the time interval $\left( 0,T\right)$ as
  described in section \ref{sec:fem}. Start with the initial approximation
  $\varepsilon_{h}^{0}= \varepsilon_g$ and compute the sequences of
  $\varepsilon_{h}^{m}$ via the following steps:

\item Step 1.  Compute solutions $E_{h}\left(
  x,t,\varepsilon_{h}^{m}\right) $ and $\lambda _{h}\left(
  x,t,\varepsilon_{h}^{m}\right) $ of state  (\ref{model1})
  and adjoint  (\ref{adjoint1}) 
problems on $K_{h}$ and $J_{\tau}$.

\item Step 2.  Update the coefficient
  $\varepsilon_h:=\varepsilon_{h}^{m+1}$  on
  $K_{h}$ and $J_{\tau}$ using the conjugate gradient method
\begin{equation}\label{cgm}
\begin{split}
\varepsilon_h^{m+1} &=  \varepsilon_h^{m}  + \alpha^m d^m(x),
\end{split}
\end{equation}
where
\begin{equation*}
\begin{split}
 d^m(x)&=  -g^m(x)  + \beta^m  d^{m-1}(x),
\end{split}
\end{equation*}
with
\begin{equation*}  
\begin{split}
 \beta^m &= \frac{\| g^m(x)\|^2}{\| g^{m-1}(x)\|^2},
\end{split}
\end{equation*}
where $d^0(x)= -g^0(x)$. In (\ref{cgm})
 the step size $\alpha$  in the gradient update is  computed as
\begin{equation}
\alpha^m = -\frac{((g^m, d^m)) }{\gamma^m {\norm{d^m}}^2},
\end{equation}
and the regularization parameter $\gamma$ is computed iteratively accordingly to \cite{BKS} as
\begin{equation}\label{iterreg}
\gamma^m = \frac{\gamma_0 }{ (m+1)^p}, p \in (0,1).
\end{equation}

\item Step 3. Stop computing $\varepsilon_{h}^{m}$ and obtain the
  function $\varepsilon_h$ at $M=m$ if either $\| g^{m}\|_{L_{2}( D_{FEM})}\leq
  \theta$ or norms $\|g^{m}\|_{L_{2}(D_{FEM})}$ are stabilized. Here
  $\theta$ is the tolerance in updates $m$ of gradient
  method. Otherwise set $m:=m+1$ and go to step 1.

\item Step 4.  Refine the mesh $K_h$ where 
\begin{equation}
|  g^M(x) | \geq  C \max_{x \in {D}_{FEM}} |  g^M(x)|,  
\end{equation}
where the constant  $C \in \left( 0,\,1\right) $ is chosen by
the user.

\item Step 5. Construct a new mesh $K_{h}$ in $D_{FEM}$ and a new partition $
J_{\tau}$ of the time interval $\left( 0,\,T\right) $. On $J_{\tau}$ the new time
step $\tau $ should be chosen in such a way that the CFL condition is
satisfied. 

\item Step 6. Interpolate the initial approximation $\varepsilon_g$
  from the previous space mesh to the new one. Set $m=1$ and return to
  step 1. 

\item Step 7. Stop  refinements of $K_h$ if norms defined in step 3 either increase or
  stabilize, compared with the previous space mesh.

\end{itemize}

\textbf{Remark}

In our computations at step 4 of the adaptive algorithm we refine only
a such domain of $D_{FEM}$ which should be designed since we assume
that we know in advance the dielectric permittivity in all other parts of $D_{FEM}$.

\section{Numerical Studies}
\label{sec:Numer-Simul}

\begin{figure}[tbp]
 \begin{center}
 \begin{tabular}{c}
{\includegraphics[scale=0.5, clip=]{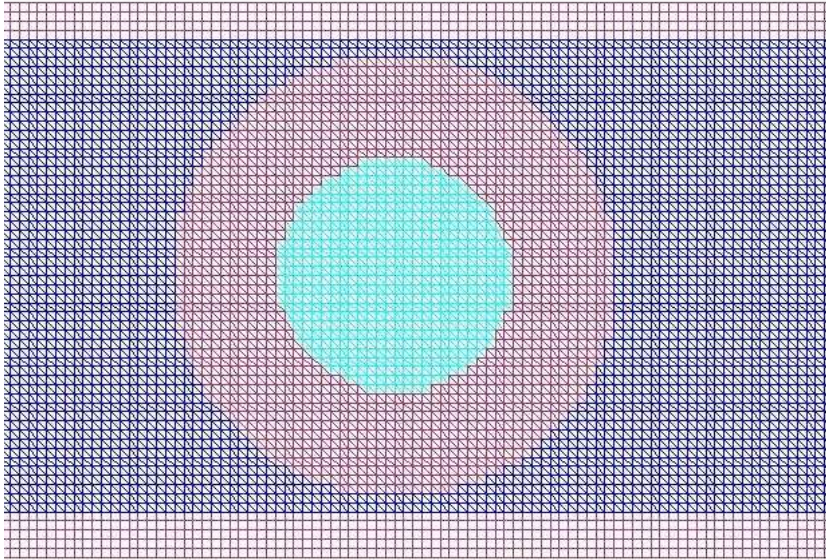}} \\
a) $D$ \\
{\includegraphics[scale=0.3, clip=]{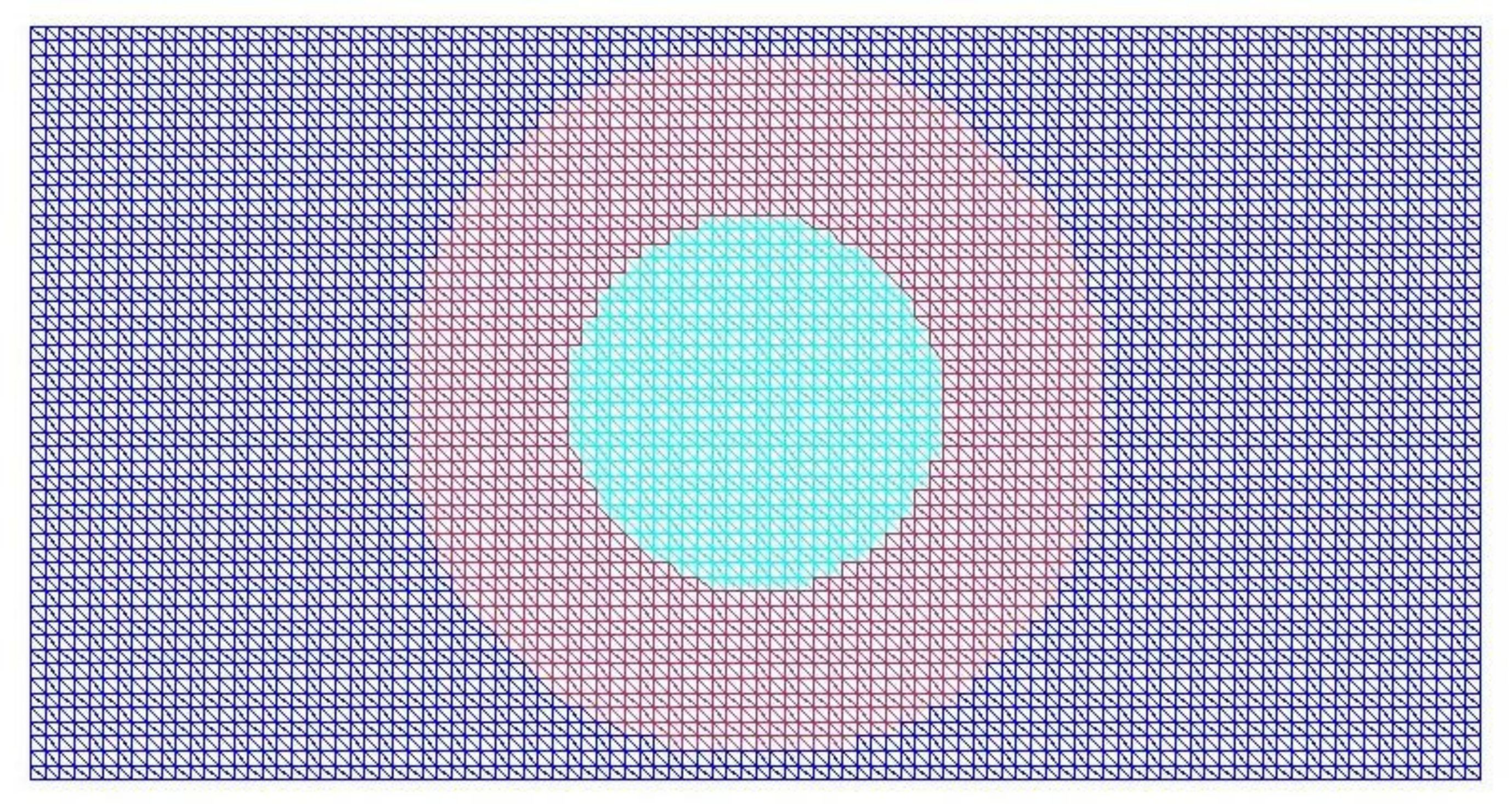}}\\
  b)   $D_{FEM}$
 \end{tabular}
 \end{center}
 \caption{ \emph{Computational coarse FE/FD mesh  used in the domain
     decomposition in $D$. b) The finite element
     mesh in $D_{FEM}$.}}
 \label{fig:0_1}
 \end{figure}

The goal of this section is to present possibility of the
computational design of nanophotonic structures with some prescribed
property.  We have chosen to design a  structure which have a
property to generate a small reflections as possible.  This problem is
equivalent to \textbf{IP}. Thus, we will reconstruct a function
$\varepsilon(x)$ inside a domain $D_{FEM}$ using the ACGM algorithm of
section \ref{subsec:ad_alg}.  We assume, that this function is known
inside $D_{FDM}$ and is set to be $\varepsilon(x)=1$.  Moreover, we
decompose also the domain $\D_{FEM}$ into three different domains
$D_1, D_2, D_3$ such that $\D_{FEM} = D_1 \cup D_2 \cup D_3$ which are
intersecting only by their boundaries, see Figure \ref{fig:0_1}.  The
boundary of $D_{FEM}$ we define as $\partial D_{FEM}$, and the boundary of
$D_1$ we define as $\partial D_1$.
  The goal of our
numerical tests is to reconstruct the dielectric permittivity function
of the approximately cyclic domain $D_{2}$ of Figure \ref{fig:0_1} which produce a
small reflections as possible.

 In our studies we initialize a plane wave $p(t)$ as the boundary
 condition on $S_{T_1}$, see (\ref{f}). Initial conditions in
 (\ref{model1}) are set to be zero.  In all computations we used the
 domain decomposition method of \cite{hybrid} implemented in the
 software package WavES \cite{waves}.  Our computational geometry $D$
 is split into two geometries $D_{FEM}$ and $D_{FDM}$ as described in
 section \ref{sec:modelhyb}, such that $D = D_{FEM} \cup D_{FDM}$, see
 Figure \ref{fig:0_1}.  We set the dimensionless computational domain
 $D$ as
 \begin{equation*}
 D = \left\{ x= (x_1,x_2) \in (-1.1, 1.1) \times (-0.62,0.62)\right\},
 \end{equation*}
and  the  domain $D_{FEM}$ as
 \begin{equation*}
 D_{FEM} = \left\{ x= (x_1,x_2) \in ((-1.0,1.0) \times (-0.52,0.52) \right\}.
 \end{equation*}
 The space mesh in $D_{FEM}$ and in $D_{FDM}$ consists of
 triangles and squares, respectively. 
  We choose the initial mesh size $h=0.02$
 in $D = D_{FEM} \cup D_{FDM}$, as well as in the
 overlapping regions between FE/FD  domains.

  We initialize  a plane wave $f(t)$ in the  equation (\ref{model1})
  in
 $D$ in time $T=[0,2.0]$ such that
 \begin{equation}\label{f}
 \begin{split}
 f\left( t\right) =\left\{ 
 \begin{array}{ll}
 \sin \left( \omega t \right) ,\qquad &\text{ if }t\in \left( 0,\frac{2\pi }{\omega }
 \right) , \\ 
 0,&\text{ if } t>\frac{2\pi }{\omega }.
 \end{array}
 \right. 
 \end{split}
 \end{equation}
 As the forward problem   in $D_{FDM}$ we
 solve the problem (\ref{model1}) choosing $\varepsilon=1$ and $D=D_{FDM}$, and
in $D_{FEM}$ we  solve
 \begin{equation} \label{3D_1}
 \begin{split}
 \varepsilon \frac{\partial^2 E}{\partial t^2} - \triangle E   &= 0,~ \mbox{in}~~ 
 D_{{FEM}},    \\
   E(x,0) = 0, ~~~E_t(x,0) &= 0~ \mbox{in}~~ D_{FEM},    \\
 E(x,t)|_{\partial D_{FEM}} &= E(x,t)|_{\partial D_{{FDM}_I}},\\
\partial_n E &= 0~ \mbox{on}~~ \partial D_1. 
 \end{split}
 \end{equation}
Here, $\partial D_{{FDM}_I}$ denote structured nodes of $D_{FDM}$
which have the same coordinates as nodes at $\partial D_{FEM}$, see
details in \cite{hybrid}. We note, that we use the boundary condition
$\partial_n E = 0$ on $\partial D_1$ which says that waves are not
penetrated into $D_1$.

We  also note that in $D_{FDM}$ the   adjoint problem  will be   the following wave equation  with $\varepsilon=1$ in  $D_{FDM}$:
\begin{equation}\label{adjwaveeq}
\begin{split}
\frac{\partial^2 \lambda}{\partial t^2} - \triangle \lambda   &= -  (E- \tilde{E}) z_{\delta},~ \mbox{in}~~ S_{T_1} \cup S_{T_2},    \\
  \lambda(x,T) = 0, ~~~\lambda_t(x,T) &=0~ \mbox{in}~~ D,     \\
\partial _{n} \lambda(x,t)& =0~ \mbox{on}~S_3.
\end{split}
\end{equation}
 Thus, as the adjoint problem in $D_{FDM}$ we
solve the problem (\ref{adjwaveeq}) and in $D_{FEM}$ we have to solve
\begin{equation} \label{adj3D_1}
\begin{split}
\varepsilon \frac{\partial^2 \lambda}{\partial t^2} - \triangle \lambda &= 0,~ \mbox{in}~~ D_{{FEM}},    \\
  \lambda(x,T) = 0, ~~~\lambda_t(x,T) &= 0~ \mbox{in}~~ D_{FEM},    \\
\lambda(x,t)|_{\partial D_{FEM}} &= \lambda(x,t)|_{\partial D_{{FDM}_{II}}}, \\
\partial_n \lambda &= 0, ~ \mbox{on}~~ \partial D_1.
\end{split}
\end{equation}
Here,  $\partial D_{{FDM}_{II}}$ denote the inner boundary of $D_{FDM}$, see details in \cite{hybrid}.

As initial guess $\varepsilon_g(x)$ we take different constant values
of the function $\varepsilon(x)$ inside domain of $D_2$ of Figure
\ref{fig:0_1} on the coarse non-refined mesh, and we take
$\varepsilon(x)=1.0$ everywhere else in $D$. We choose three different
constant values of $\varepsilon_g(x)= \{ 0.5,1.5,2.0,2.5\}$ inside $D_2$.
  We define that the minimal and maximal values of the
   function $\varepsilon(x)$  belongs
   to the following set $M_{\varepsilon} $ of admissible parameters
 \begin{equation}\label{admpar}
 \begin{split}
  M_{\varepsilon} \in \left \{\varepsilon\in C(\overline{D })| \frac{1}{\max_{D_2} \varepsilon_g(x)} \leq \varepsilon(x)\leq  \max_{D_2} \varepsilon_g(x) \right \}.
 \end{split}
 \end{equation}
 The time step is chosen to be $\tau=0.002$ which satisfies
the CFL condition \cite{CFL67}.

 \subsection{Reconstructions}

 We generate data at the observation points at $S_{T_1} \cup S_{T_2}$
 by solving the forward problem (\ref{model1}) in the time interval
 $t=[0,2.0]$, with function $f(t)$ given by (\ref{f}) and $\omega=40$.
 To generate $\tilde{E}$ at $S_{T_1} \cup S_{T_2}$ we take the
 function $\varepsilon(x)=1$ for all $x$ in $D$ and solve the problem
 (\ref{model1}) with a plane wave (\ref{f}) and $\omega=40$.

We regularize the solution of the inverse problem by starting
computations with regularization parameter $\gamma=0.01$ in
(\ref{functional}) and then updating this parameter iteratively in ACGM by
formula (\ref{iterreg}).  Computing of the regularization parameter by this way is
optimal one  for our problem.  We
refer to \cite{Engl} for different techniques for choice of a
regularization parameters.

\begin{figure}[tbp]
\begin{center}
\begin{tabular}{cc}
{\includegraphics[scale=0.2, clip=]{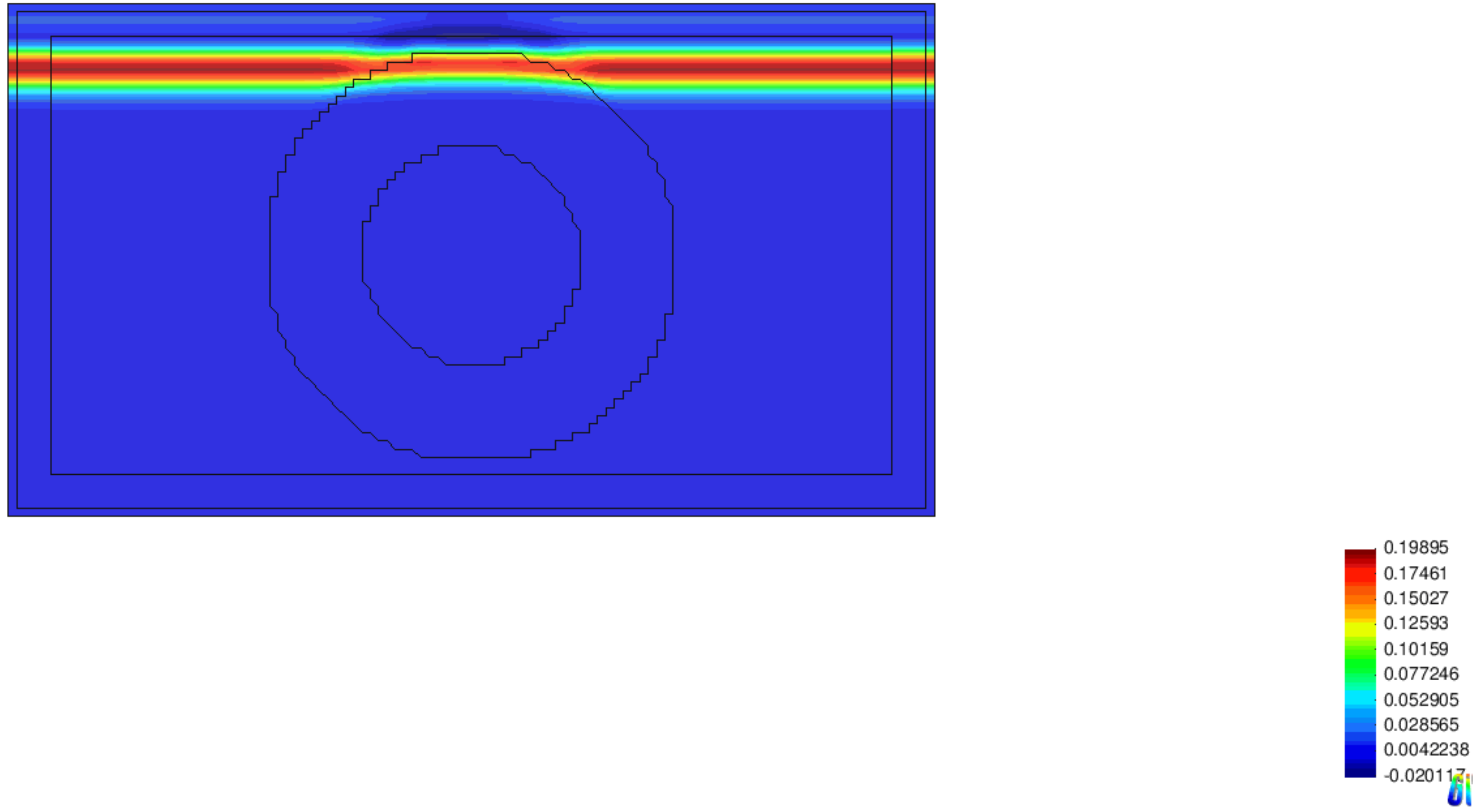}} &
{\includegraphics[scale=0.2, clip=]{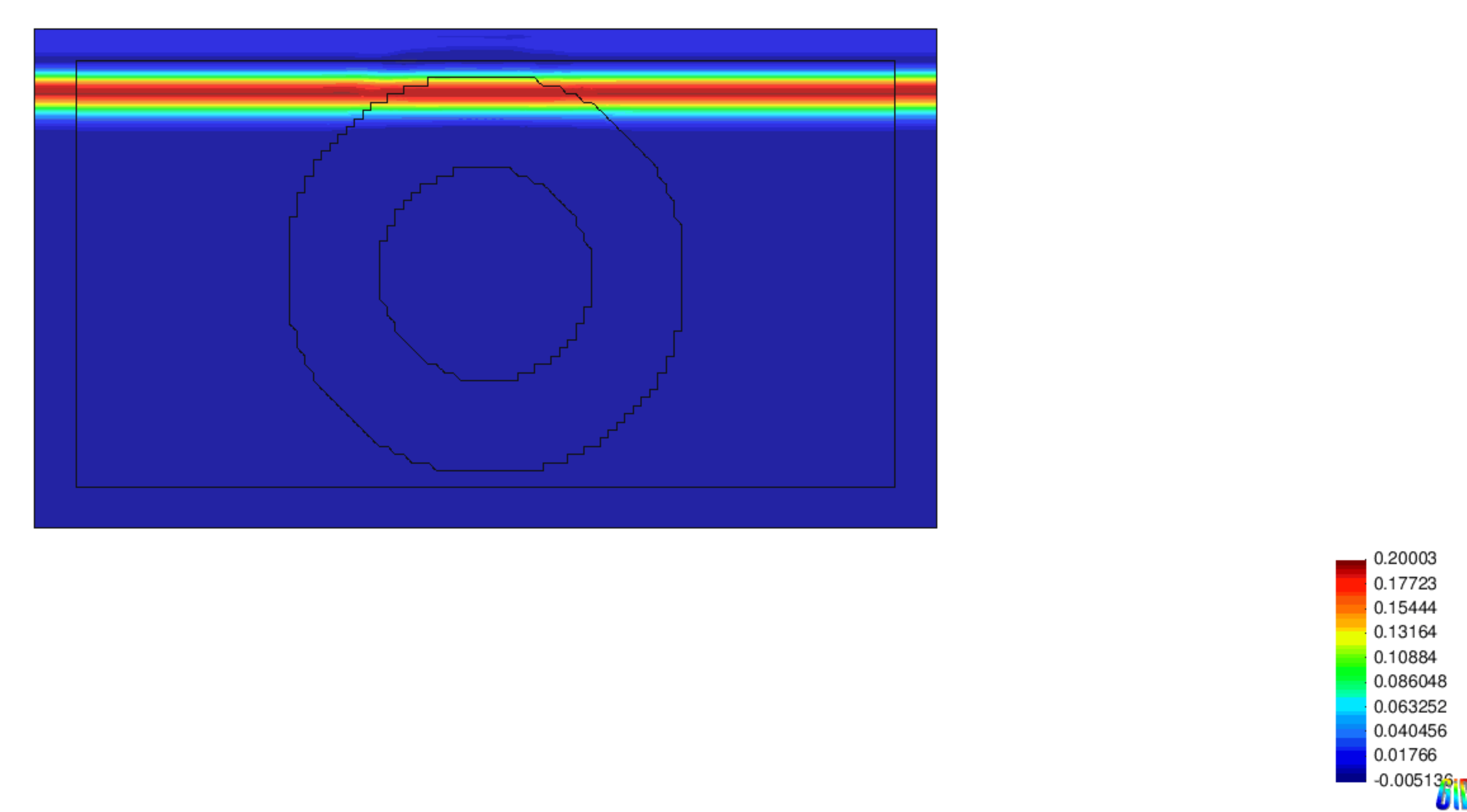}} \\
a)  t= 0.3 & b) t= 0.3\\
{\includegraphics[scale=0.2, clip=]{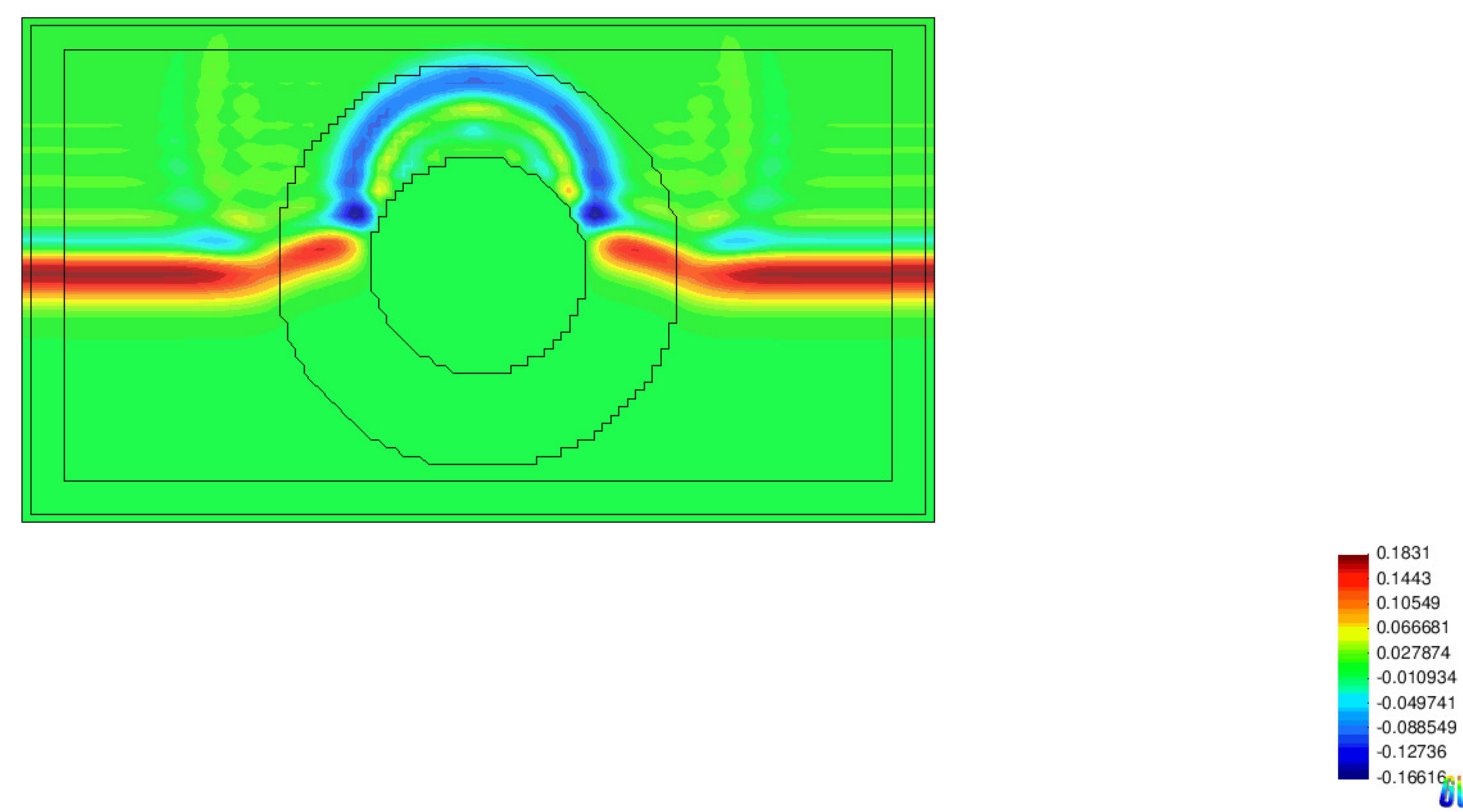}} &
{\includegraphics[scale=0.2, clip=]{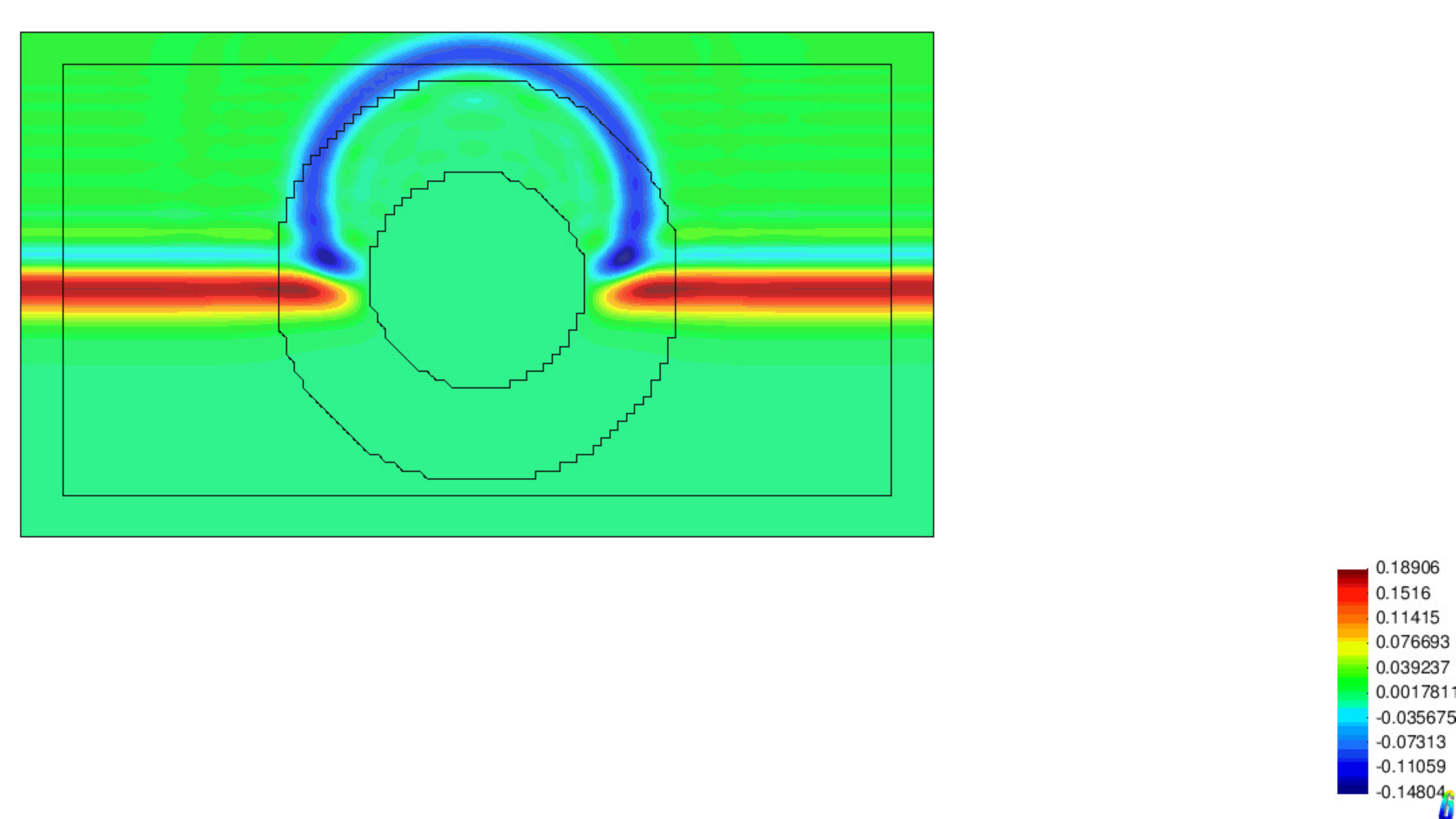}} \\
c)  t= 0.78 & d)  t= 0.78 \\
{\includegraphics[scale=0.2, clip=]{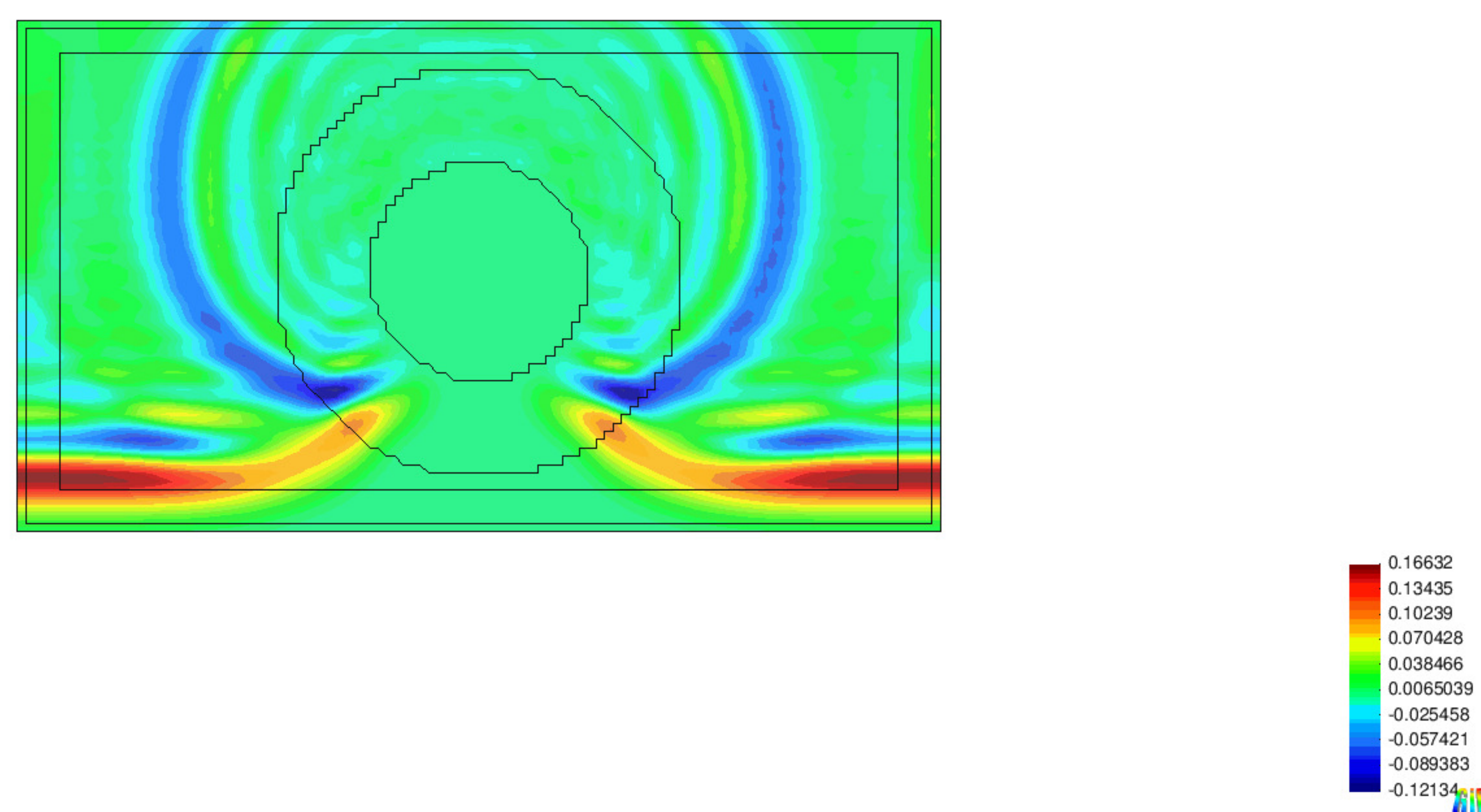}} &
{\includegraphics[scale=0.2, clip=]{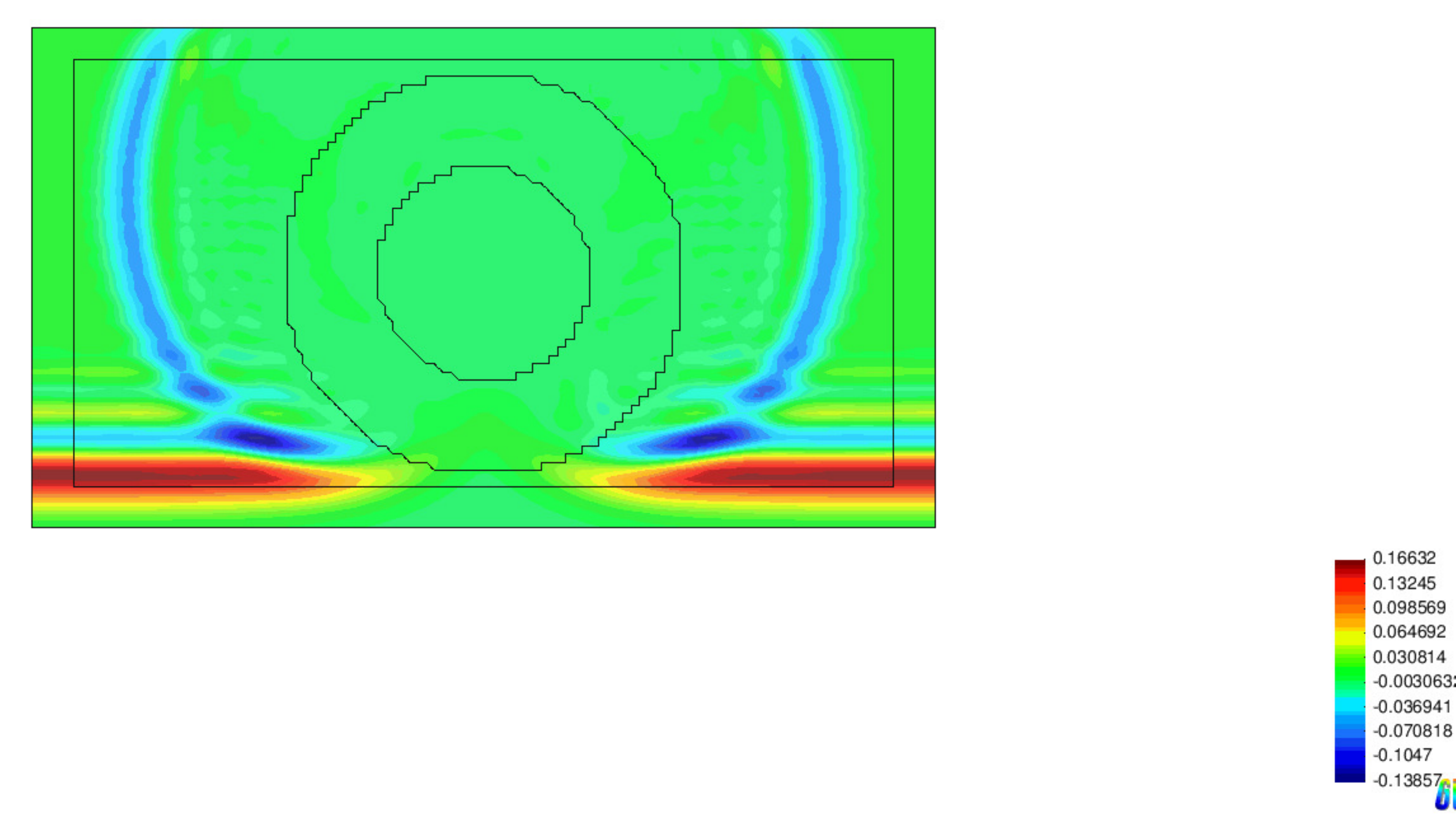}} \\
e)  t= 1.26 & f) t= 1.26\\
\end{tabular}
\end{center}
\caption{{\protect\small \emph{Computational solution of
      (\ref{model1}) using domain decomposition method of
      \cite{hybrid} at different times: a),c),e) on the coarse mesh
      with $\varepsilon_g=1.5$ in $D_2$; b),d),f) on the four times
      refined mesh with optimized $\varepsilon$ of Figure
      \ref{fig:4}-d).}}}
\label{fig:2}
\end{figure}

\begin{figure}[tbp]
\begin{center}
\begin{tabular}{cc}
{\includegraphics[scale=0.2, clip=]{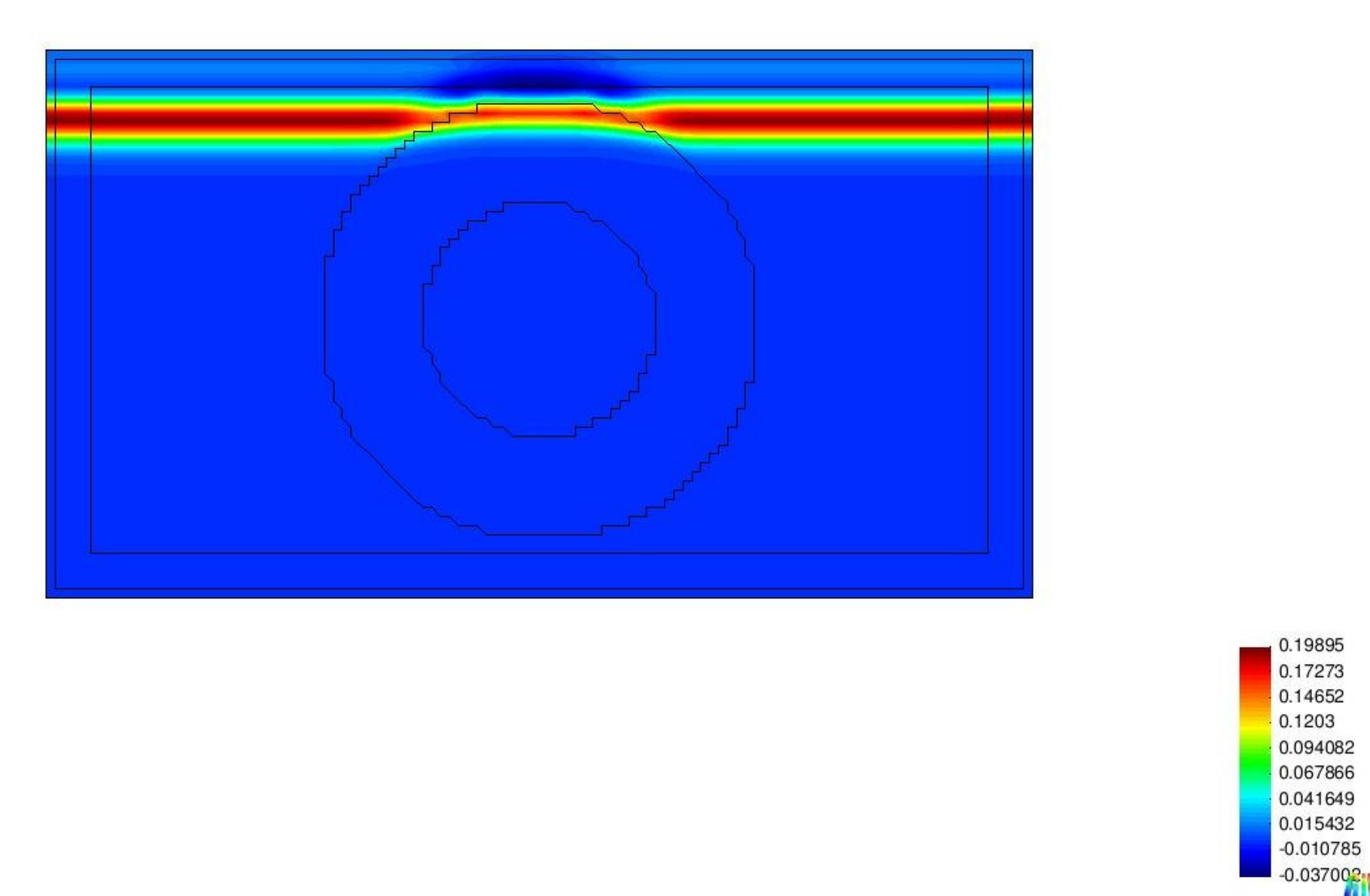}} &
{\includegraphics[scale=0.2, clip=]{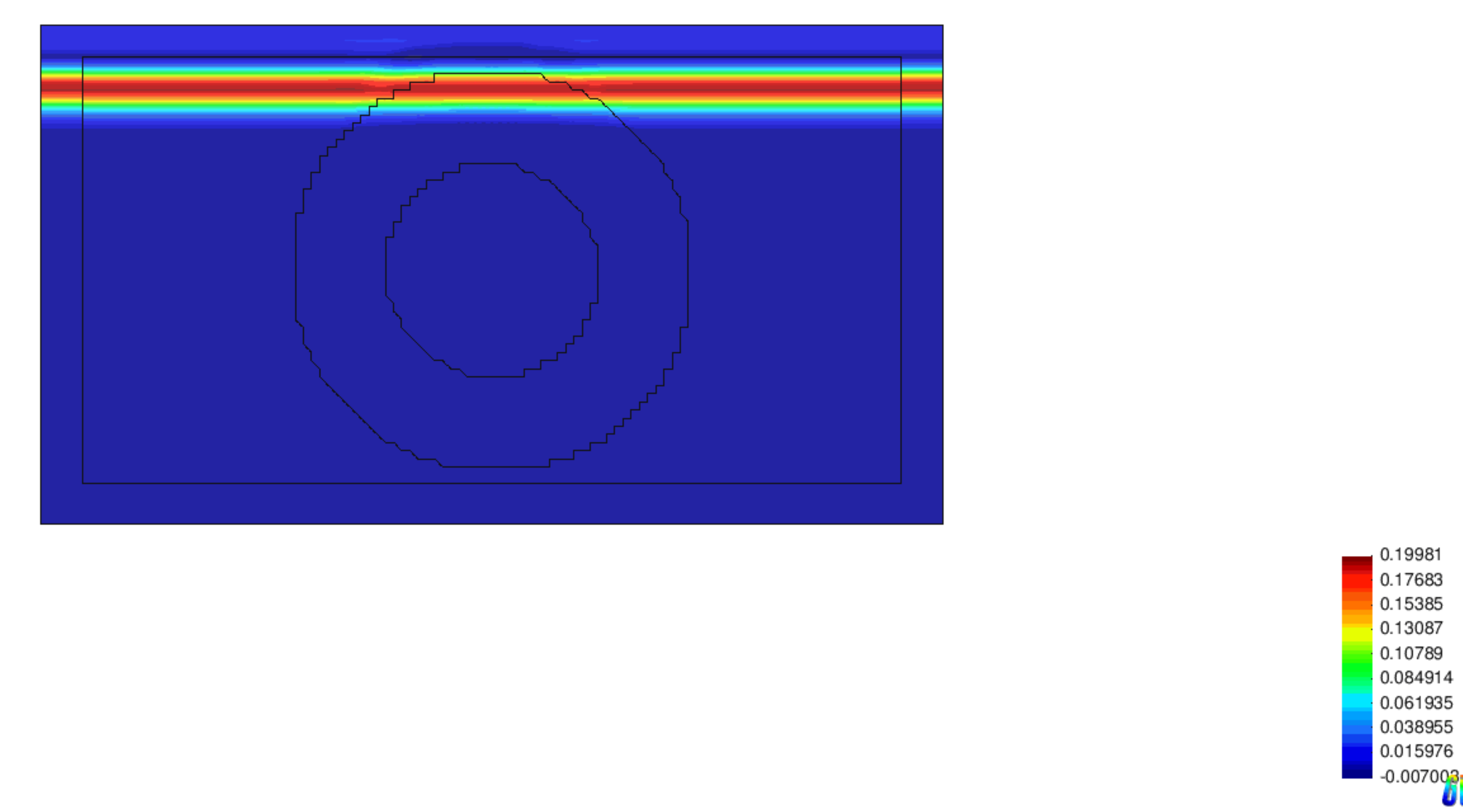}} \\
a)  t= 0.3 & b) t= 0.3\\
{\includegraphics[scale=0.2, clip=]{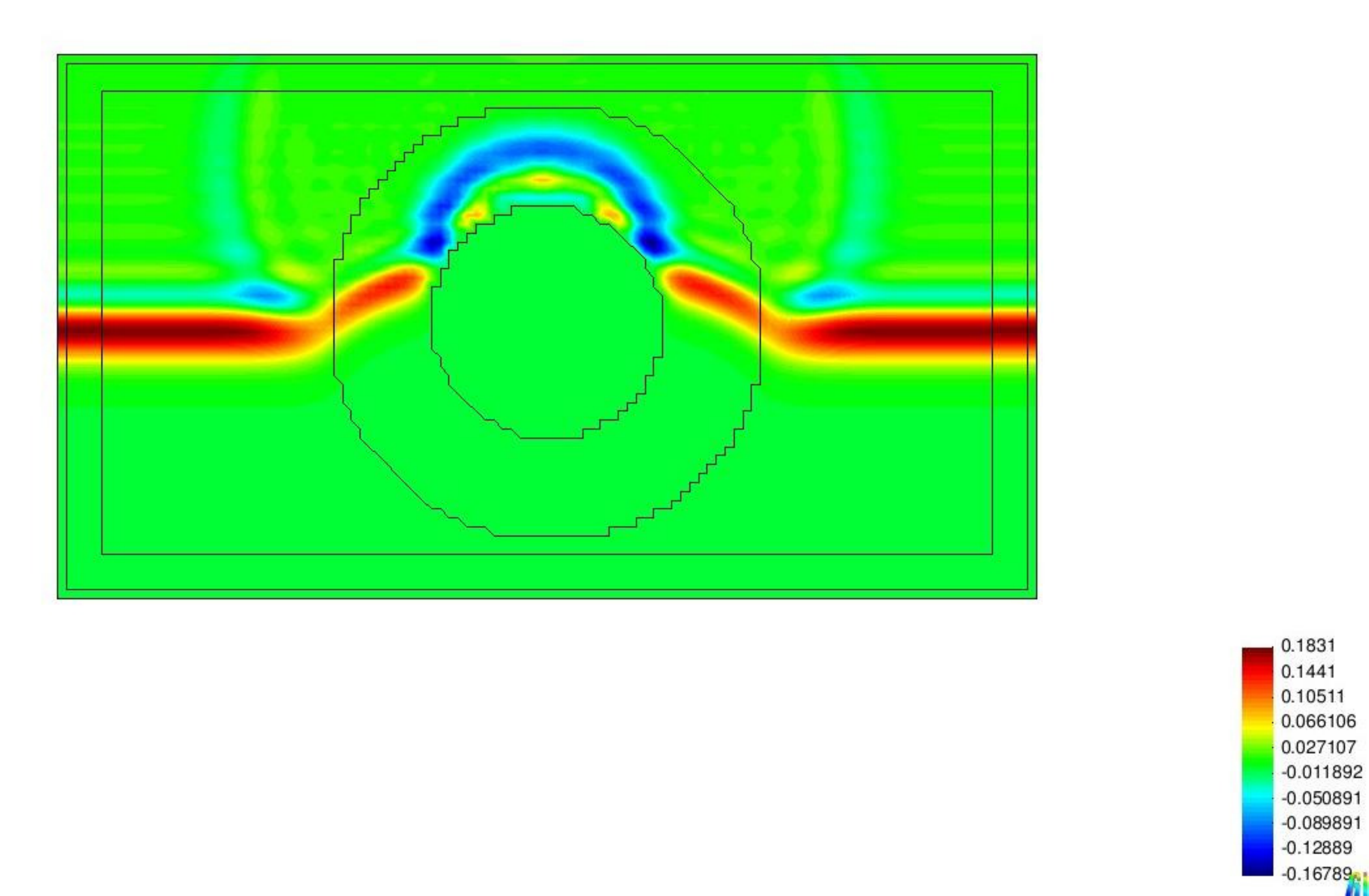}} &
{\includegraphics[scale=0.2, clip=]{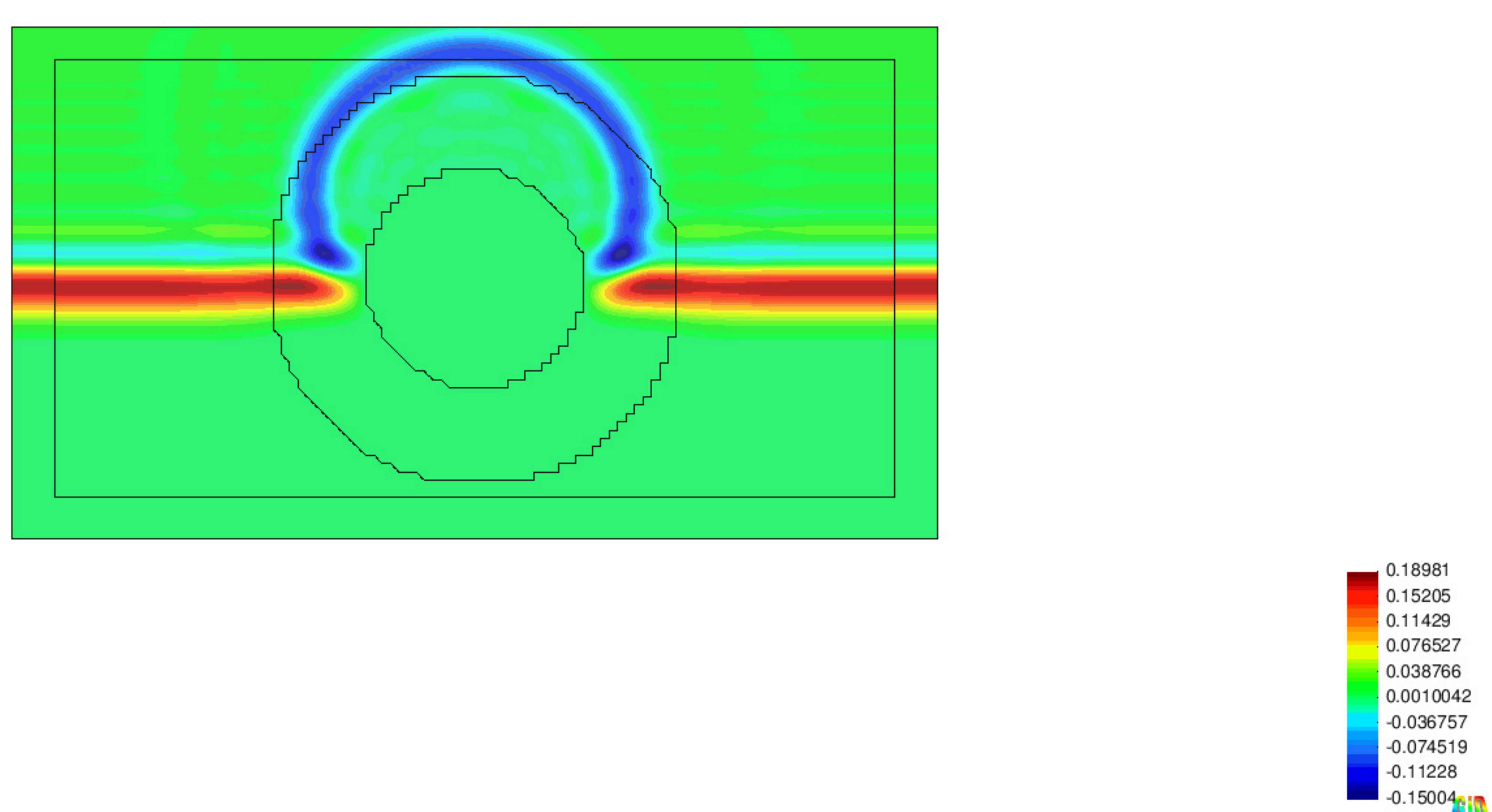}} \\
c)  t= 0.78 & d)  t= 0.78 \\
{\includegraphics[scale=0.2, clip=]{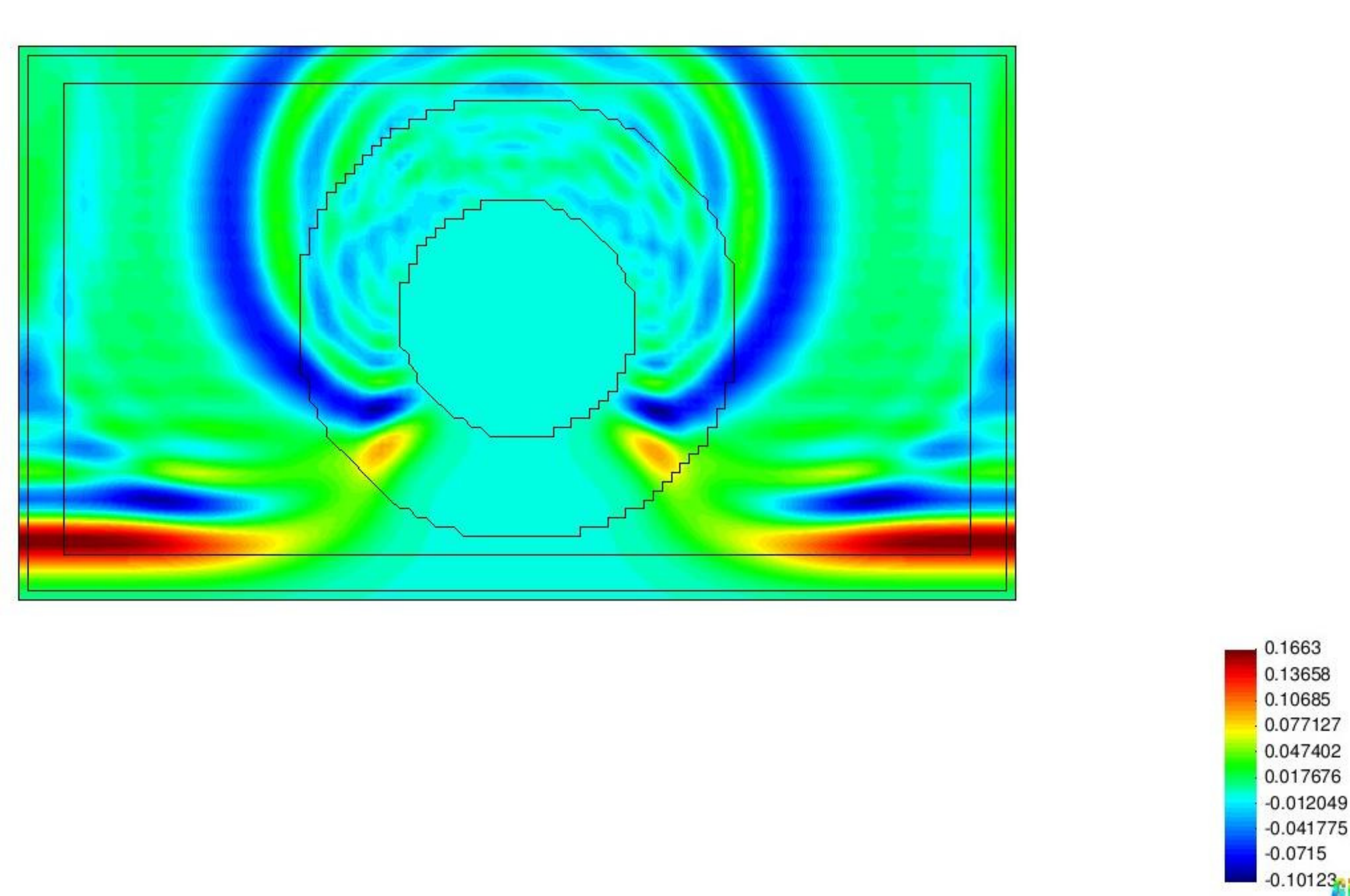}} &
{\includegraphics[scale=0.2, clip=]{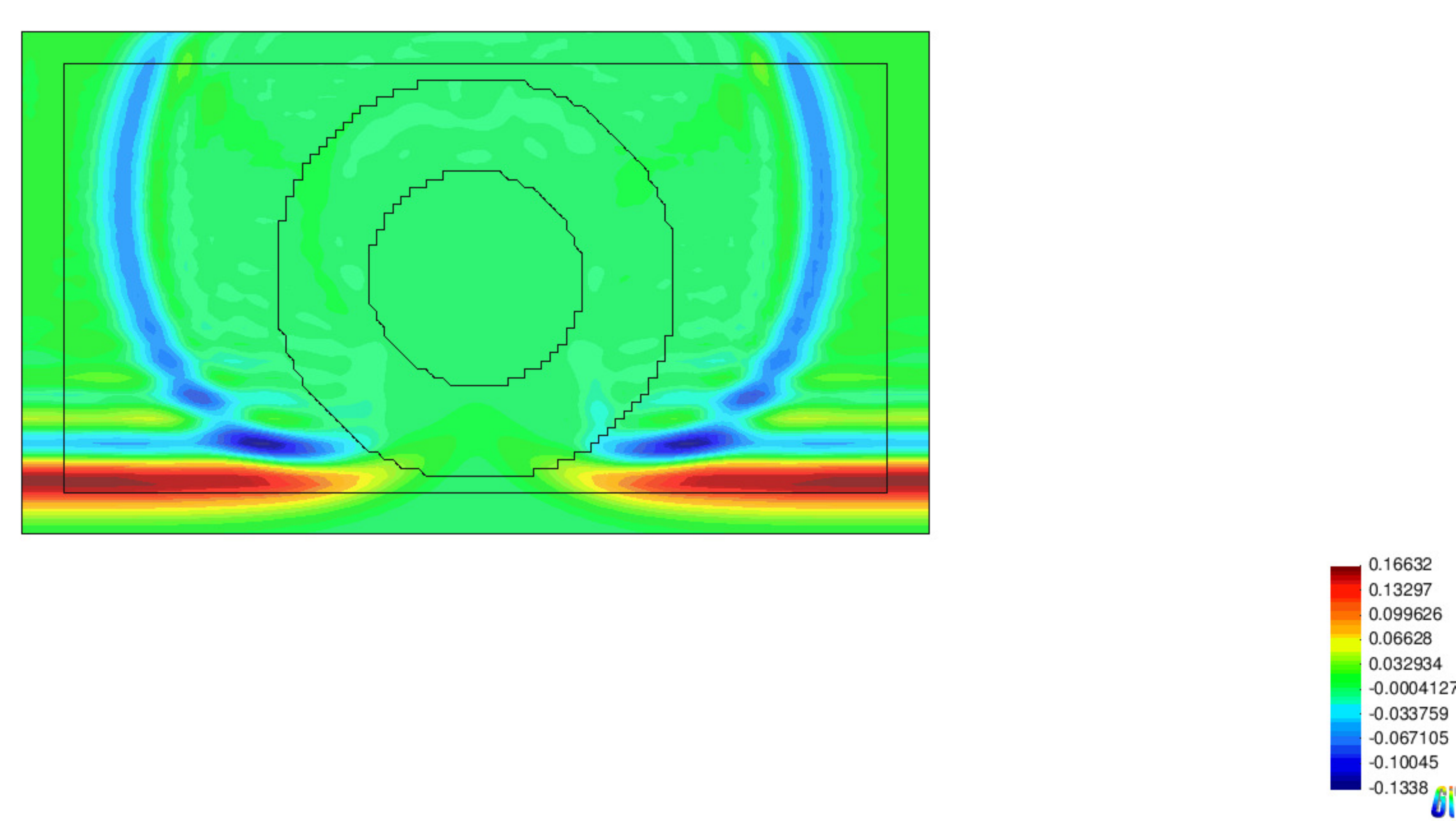}} \\
e)  t= 1.26 & f) t= 1.26\\
\end{tabular}
\end{center}
\caption{{\protect\small \emph{Computational solution of (\ref{model1})
 using domain decomposition method of \cite{hybrid}
      at different times: a),c),e) on the coarse mesh with
      $\varepsilon_g=2$ in $D_2$; b),d),f) on the four times refined mesh with optimized 
$\varepsilon$ of Figure \ref{fig:4}-f).}}}
\label{fig:3}
\end{figure}

Figures \ref{fig:2}, \ref{fig:3} show time-dependent reflections from
the dielectric permittivity function when $\varepsilon=\varepsilon_g$
(on the left) and after optimization procedure after four
refinements of the mesh in $D_2$ (on the right).  All right
figures of Figures \ref{fig:2}, \ref{fig:3} show significant reduction of reflections compared with left
figures.

\begin{figure}[tbp]
\begin{center}
\begin{tabular}{cc}
{\includegraphics[scale=0.2, clip=]{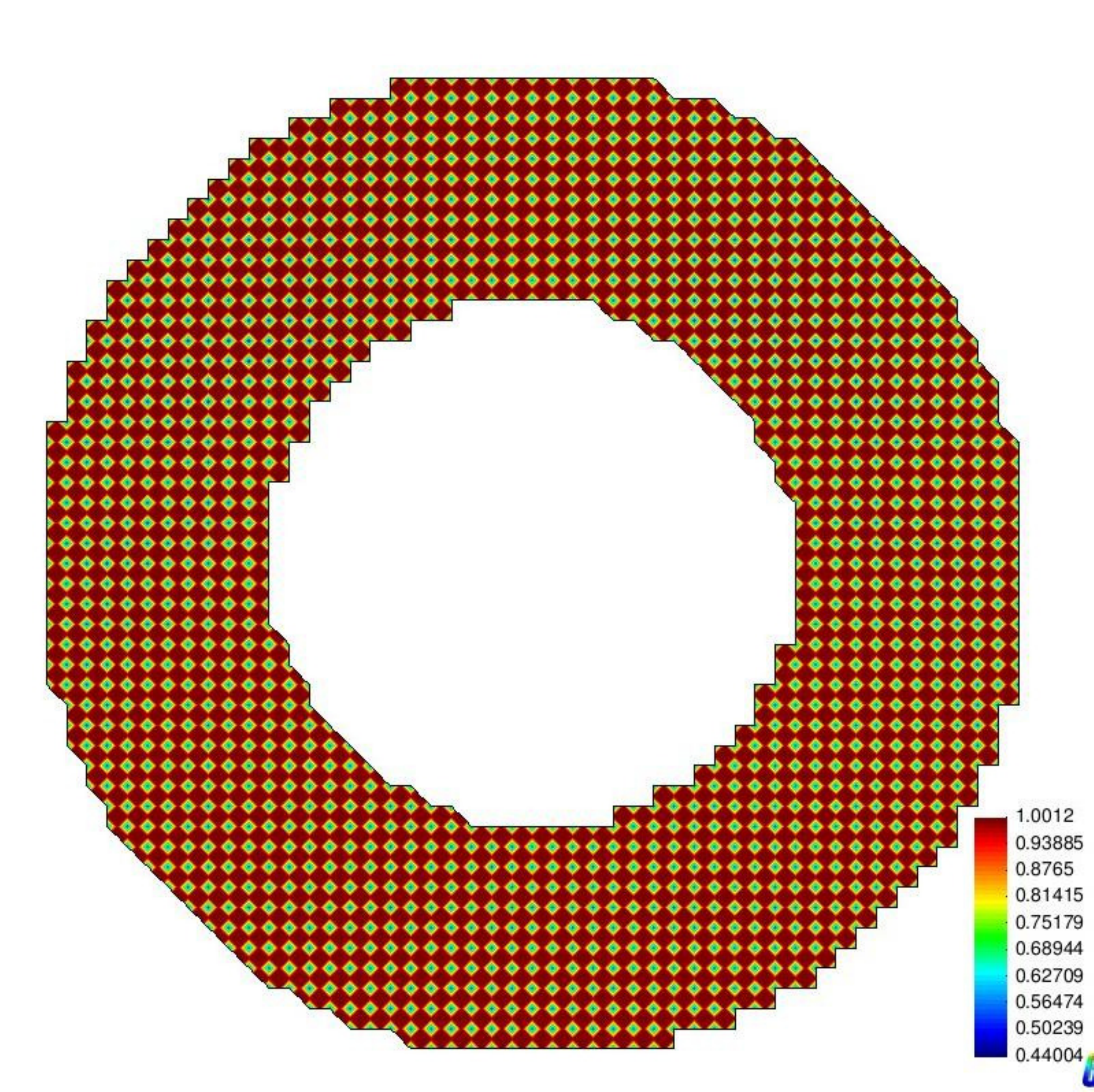}} &
{\includegraphics[scale=0.2, clip=]{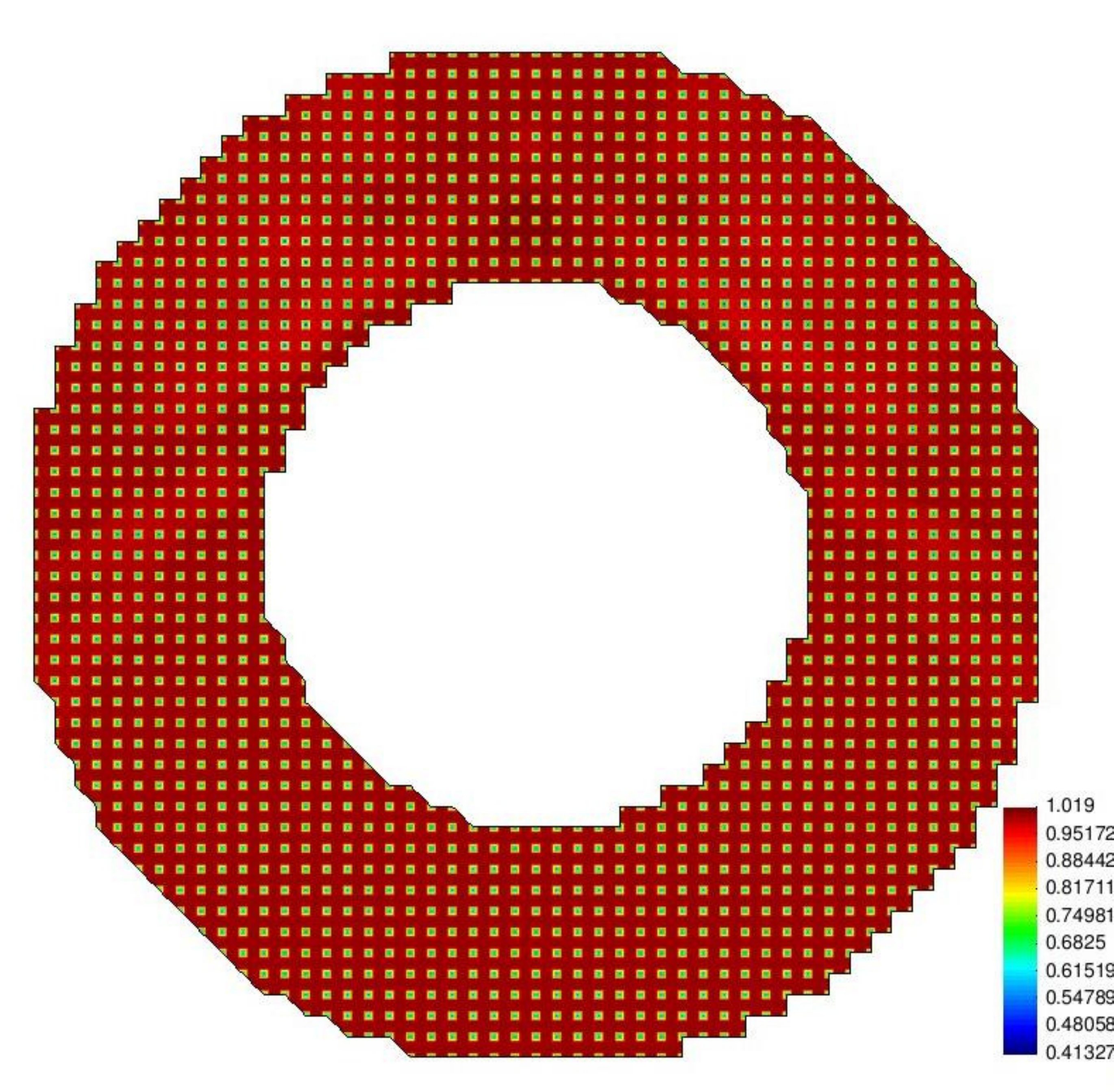}} \\
 a)  $\varepsilon_g = 0.5, n=3$  & { b)  $\varepsilon_g = 0.5, n=4$ } \\
{\includegraphics[scale=0.2, clip=]{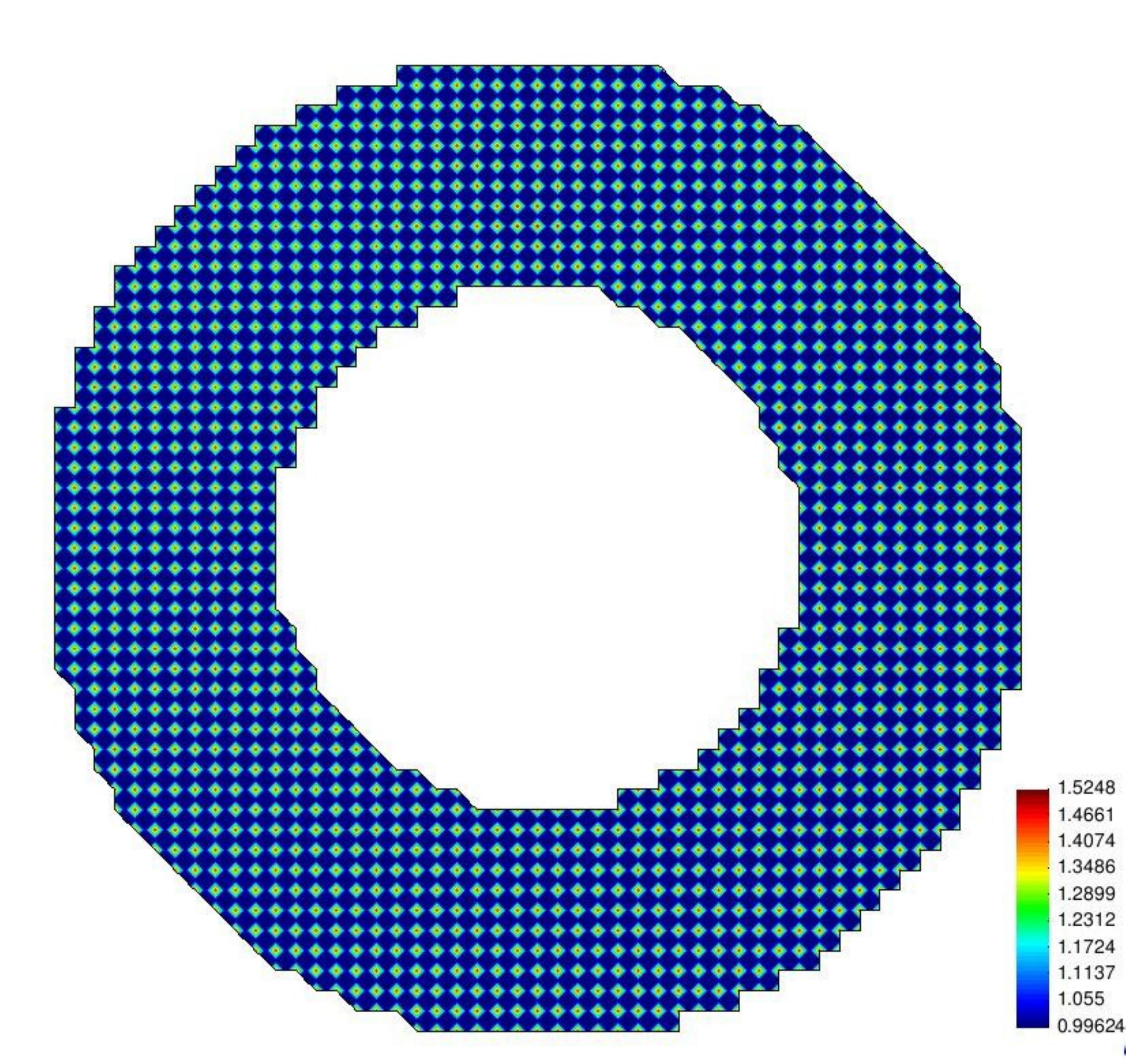}} &
{\includegraphics[scale=0.2, clip=]{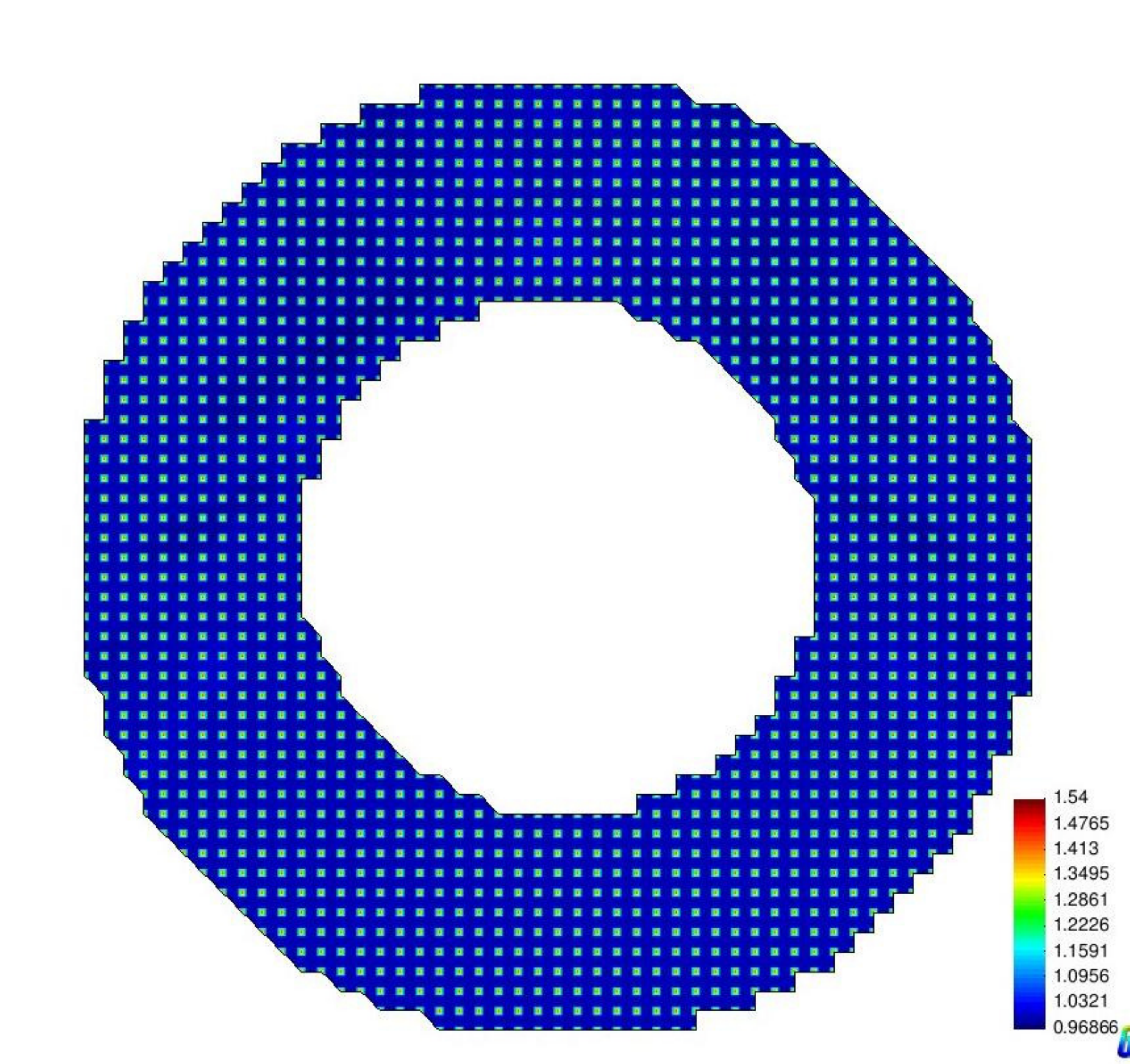}} \\
{ c)    $\varepsilon_g = 1.5, n=3$  } &{ d)   $\varepsilon_g = 1.5, n=4$   } \\
{\includegraphics[scale=0.2, clip=]{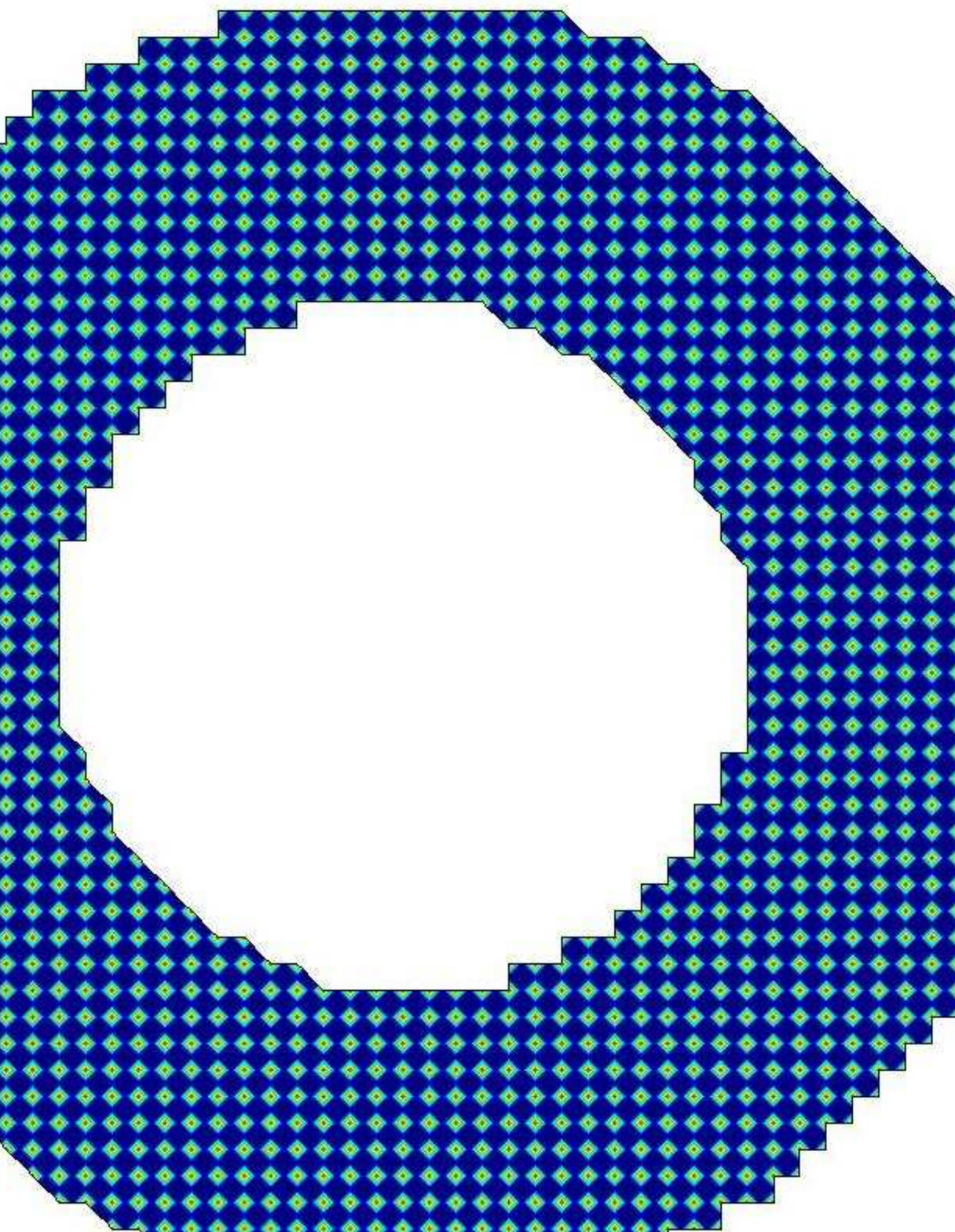}} &
{\includegraphics[scale=0.2, clip=]{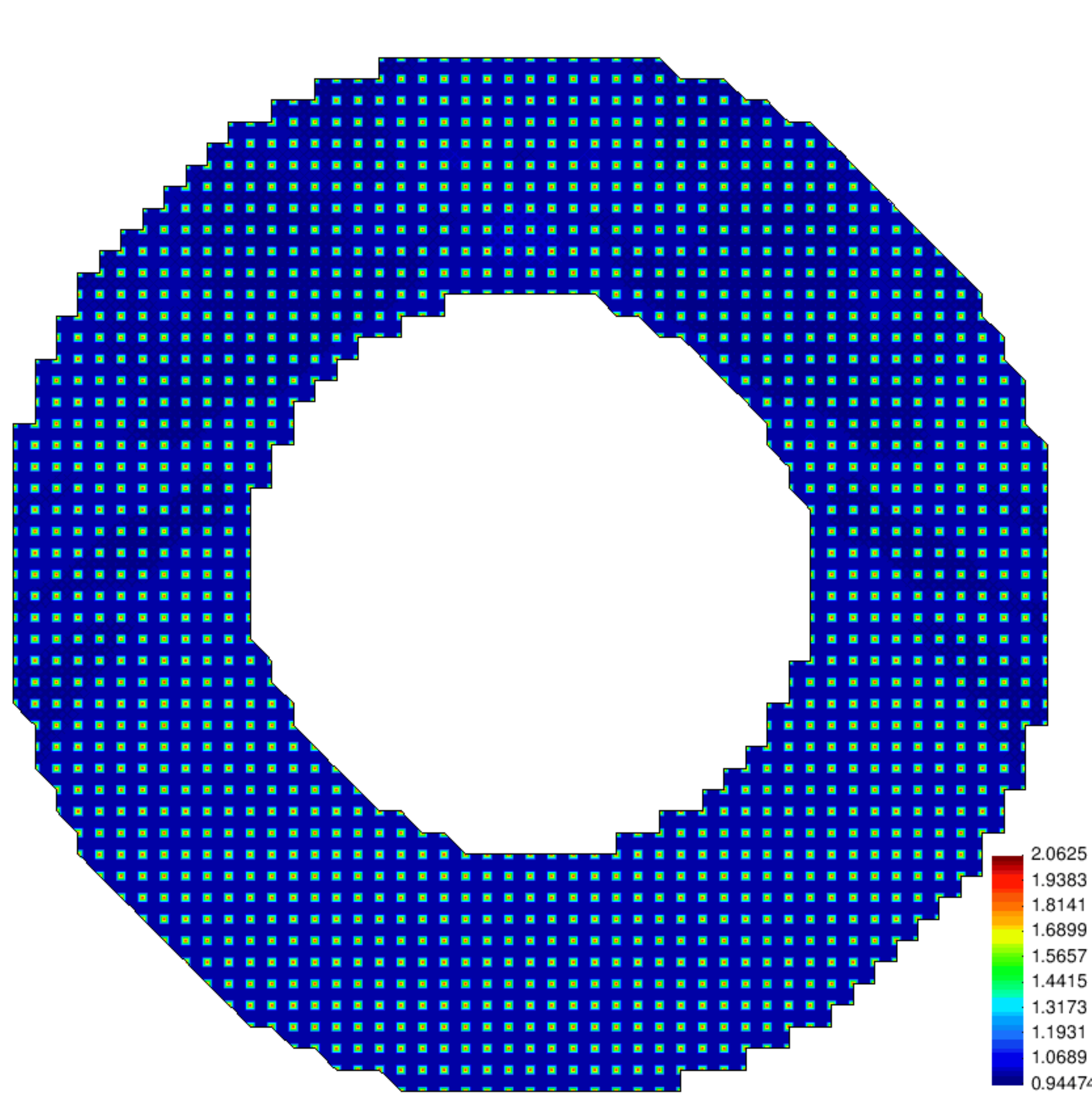}} \\
{ e)    $\varepsilon_g = 2, n=3$  } &{ f)   $\varepsilon_g = 2.0, n=4$   } 
\end{tabular}
\end{center}
\caption{{\protect\small \emph{Reconstructions in $D_2$ on three and
      four-times adaptively refined meshes for different
      $\varepsilon_g$.}}}
\label{fig:4}
\end{figure}

\begin{figure}[tbp]
\begin{center}
\begin{tabular}{cc}
{\includegraphics[scale=0.4, clip=]{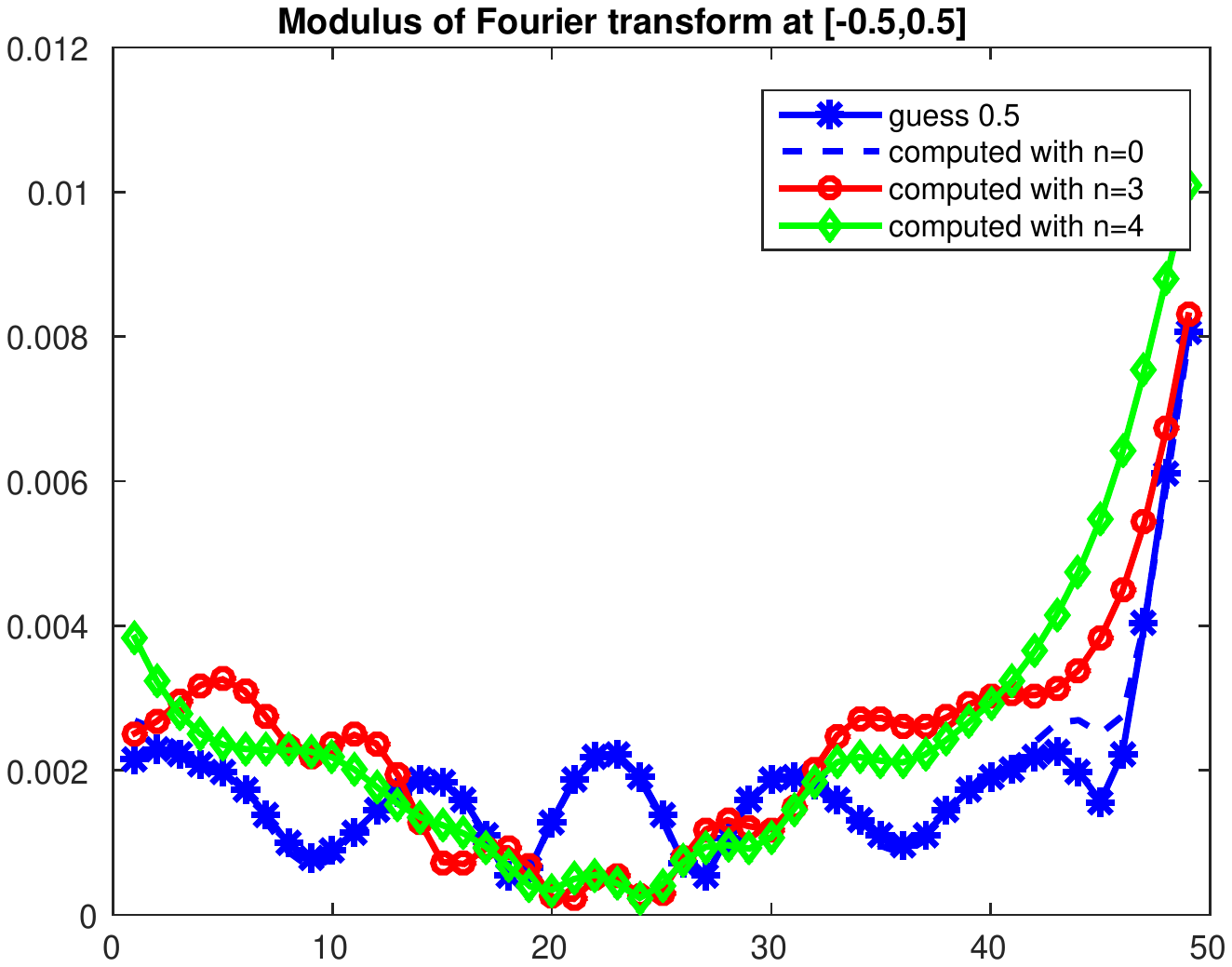}} &
{\includegraphics[scale=0.4,clip=]{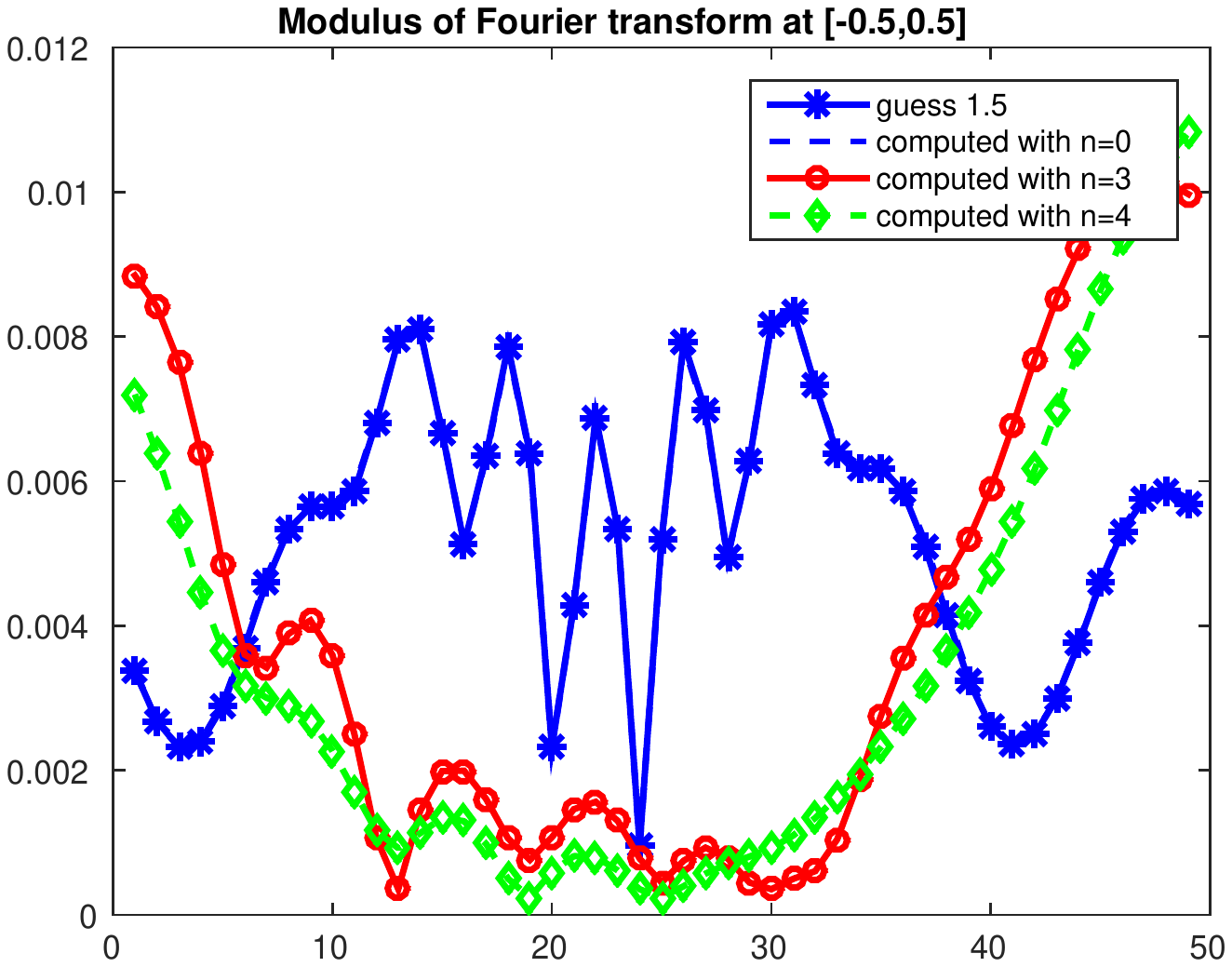}} \\
\tiny{ a)  $\varepsilon_g = 0.5$ } & \tiny{ b)   $\varepsilon_g = 1.5$ } \\
{\includegraphics[scale=0.4, clip=]{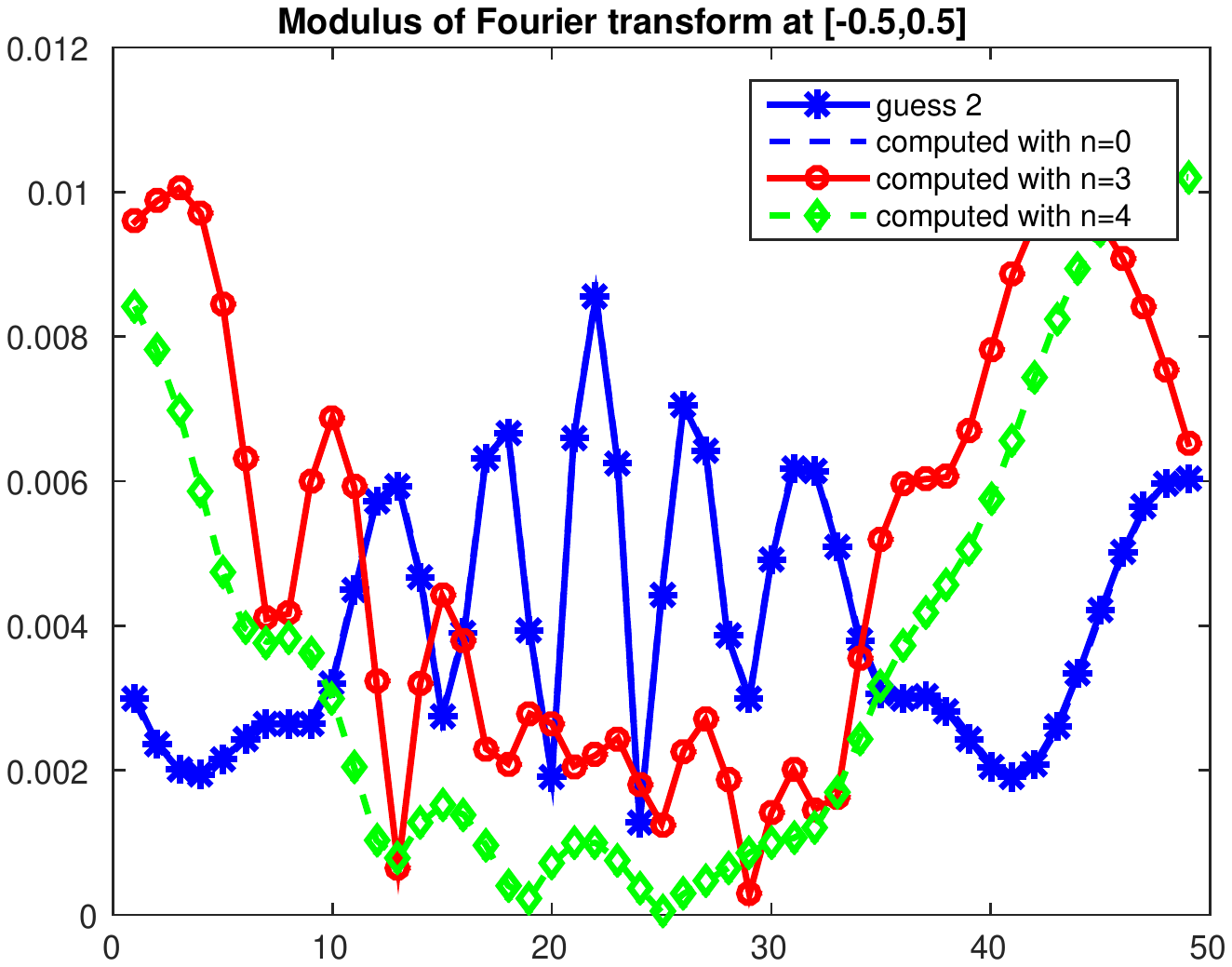}} &
{\includegraphics[scale=0.4, clip=]{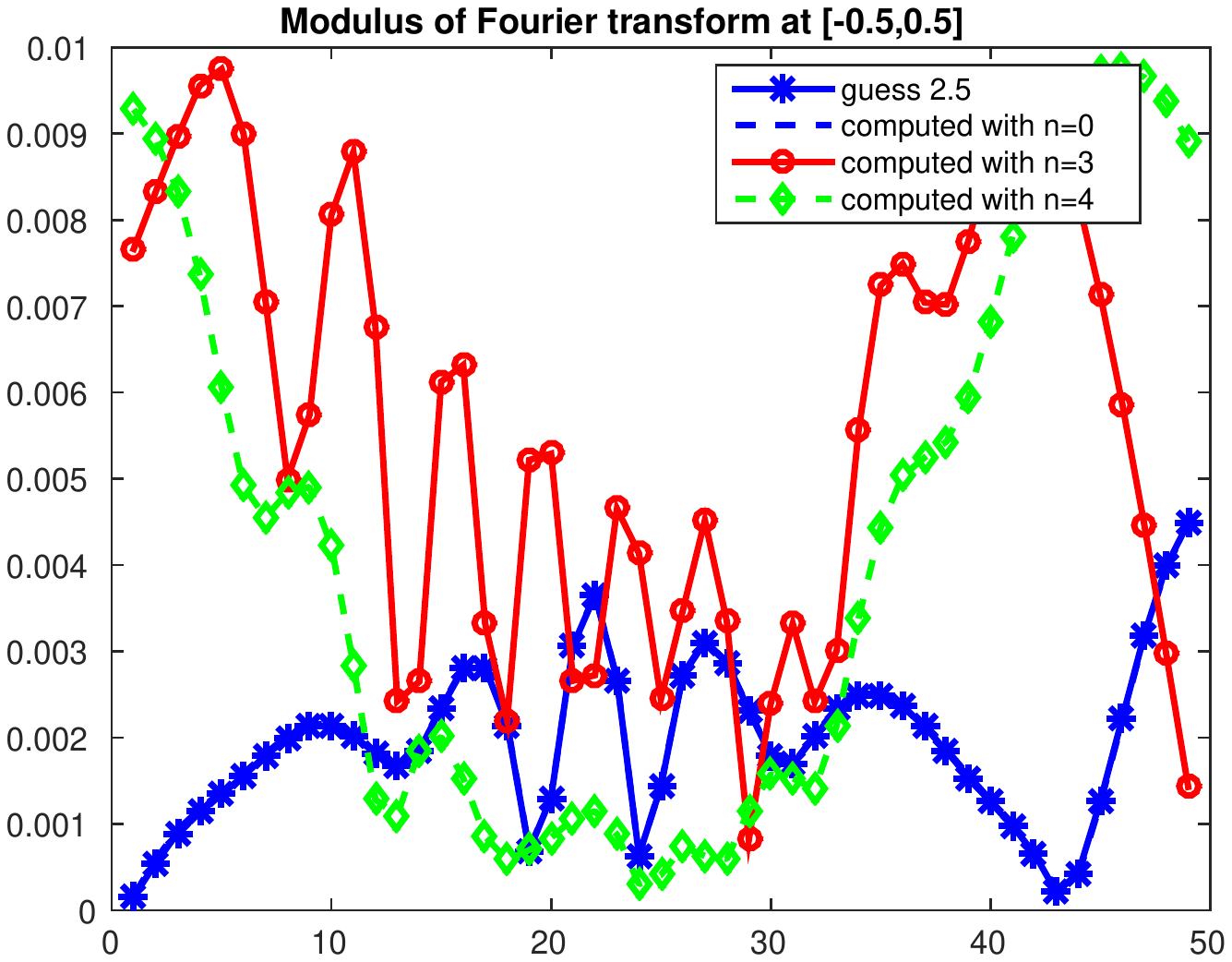}} \\
\tiny{ c)   $\varepsilon_g = 2.0 $} &\tiny{ d)   $\varepsilon_g = 2.5$} \\
\end{tabular}
\end{center}
\caption{{\protect\small \emph{Modulus of the Fourier transform at
      $\partial_2 D$ for different $\varepsilon_g$ in $D_2$ after
      applying the ACGM algorithm. Here, $n$ is the number of refinements of
      the mesh.}}}
\label{fig:5}
\end{figure}

 Figures \ref{fig:4} present reconstructions which we have obtained on
 three and four times adaptively refined mesh when we take different
 initial guesses on the coarse mesh. All guesses produce different
 structures of the domain $D_2$ with different values of the function
 $\varepsilon(x)$ inside it.  Smallest reflections we obtain taking
 the initial guess $\varepsilon_g =0.5$ inside $D_2$, and largest -
 with $\varepsilon_g =2.5$, see Figure \ref{fig:5}.
Figures \ref{fig:5} present comparison of reflections from initial and
optimized  functions $\varepsilon(x)$ after applying the Fourier
transform to the solution $E(x,t)$.

Interesting designed domains are obtained with initial guesses
$\varepsilon_g(x)= \{ 0.5,1.5,2.0\}$. In this case we obtain
optimized values of $\varepsilon(x)$ which can be of physical
interest, see Figures \ref{fig:4}-b), d), f).

 \begin{figure}[tbp]
 \begin{center}
 \begin{tabular}{c}
{\includegraphics[scale=0.25, clip=]{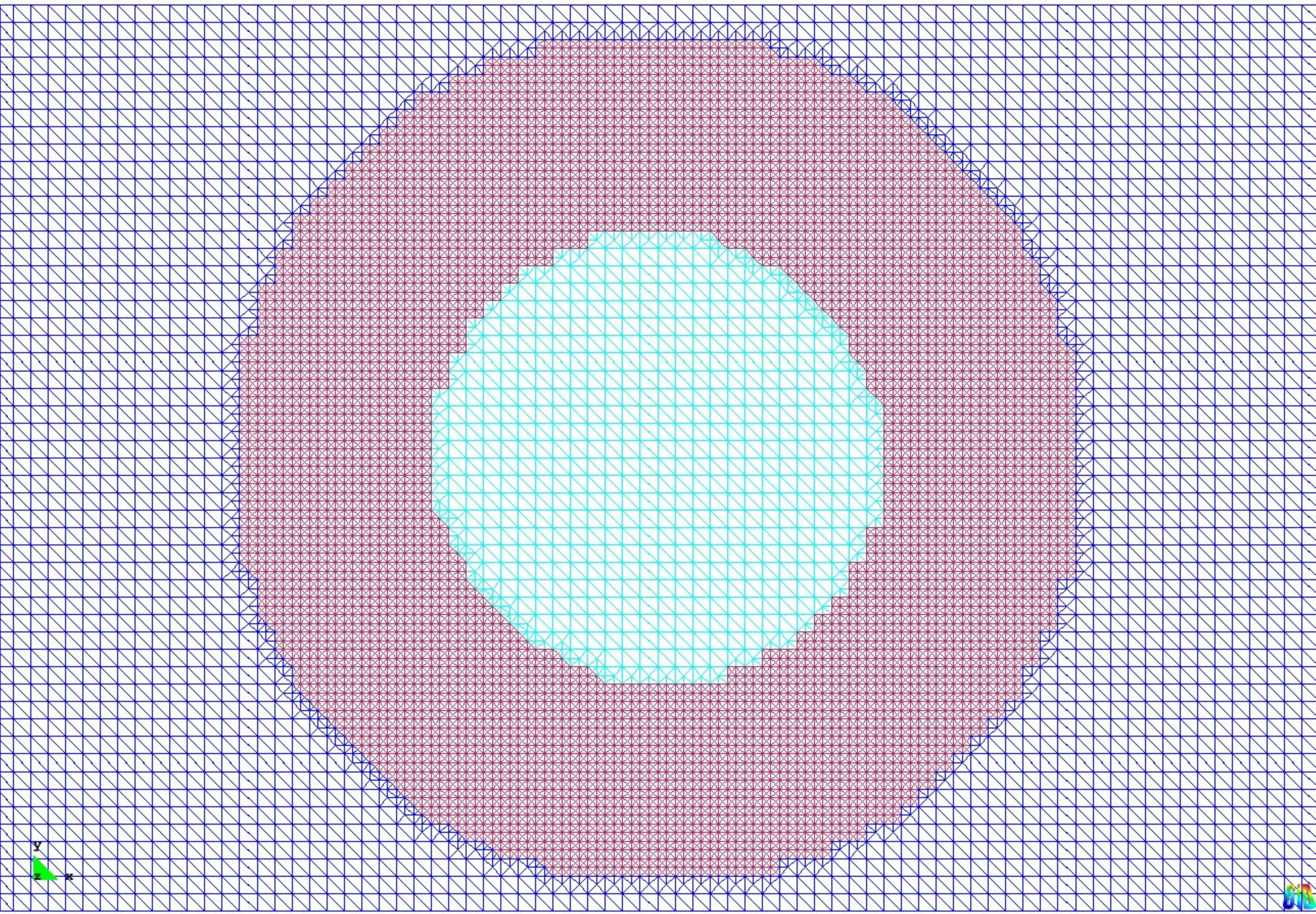}} \\
a) $n=3$ \\
{\includegraphics[scale=0.25, clip=]{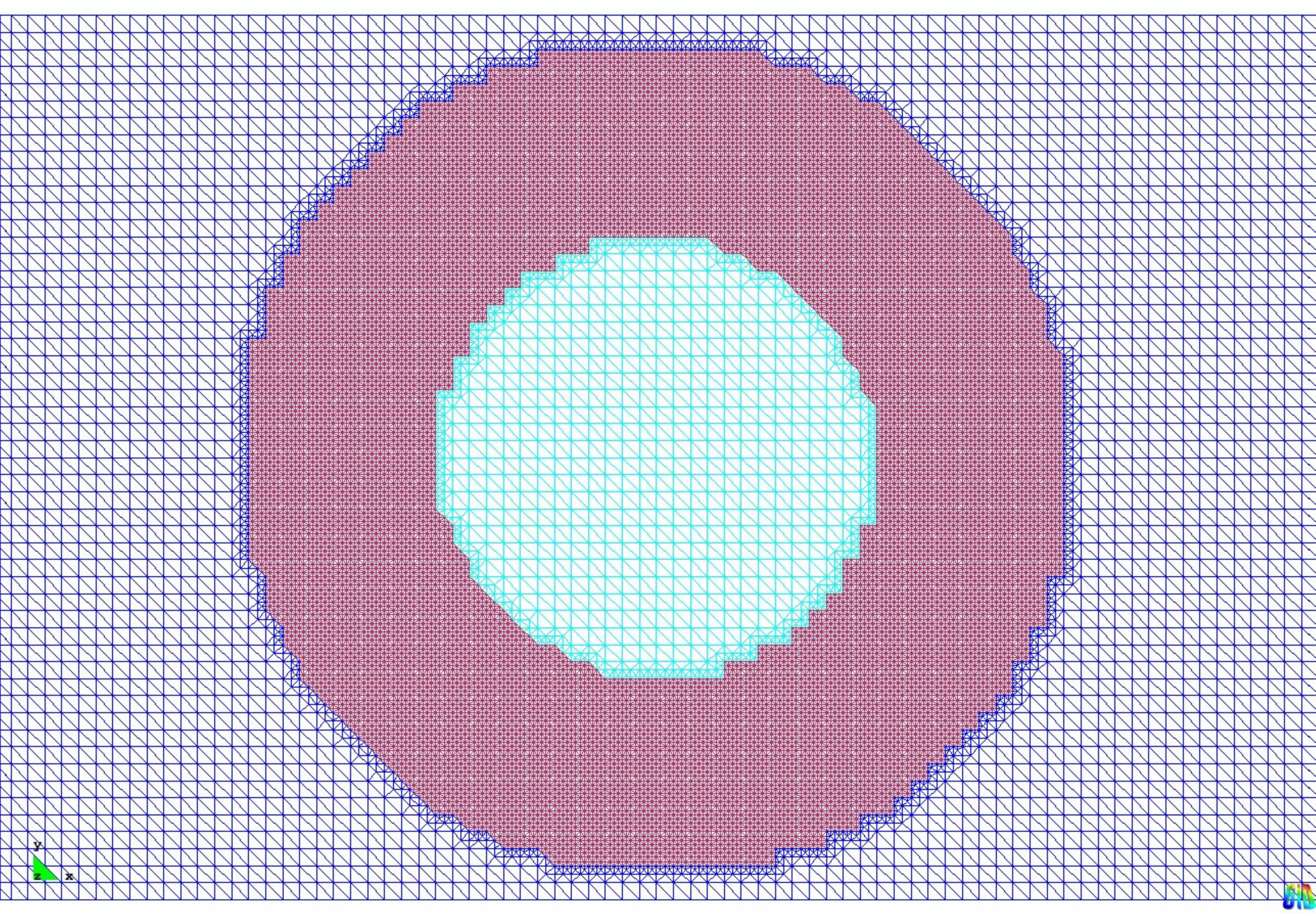}}\\
 b)   $n=4$
 \end{tabular}
 \end{center}
 \caption{ \emph{Zoomed main parts of computationally adaptively
     refined meshes for $\varepsilon_g=1.5$: a) three-times refined
     mesh; b) four-times refined mesh.}}
 \label{fig:1}
 \end{figure}

\section*{Acknowledgments}

This  work is supported by the funding from the Area of Advance
``Nanoscience and Nanotechnology''. The research of L.B is
  supported by the sabbatical programme at the Faculty of Science,
  University of Gothenburg.


\end{document}